\documentclass[12pt]{amsart}

\usepackage[dvips]{graphicx}
\usepackage{amsmath,amssymb,amsthm,wrapfig,amsfonts}

\newtheorem{defn}{Definition}[section]
\newtheorem{thm}[defn]{Theorem}
\newtheorem{prop}[defn]{Proposition}
\newtheorem{cor}[defn]{Corollary}
\newtheorem{lem}[defn]{Lemma}
\newtheorem{rem}[defn]{Remark}
\newtheorem{ex}[defn]{Example}
\newtheorem{claim}[defn]{Claim}

\newtheorem{ack}{Acknowledgements\!\!}

\DeclareMathOperator{\mint}{
\mathchoice{%
\ooalign{%
\ensuremath{%
\begin{picture}(8,8)
\thicklines
\put(1.25,3.0){\line(1,0){8}}
\end{picture}}%
\crcr
\hss\ensuremath{\displaystyle\int}\hss}
\hspace*{-4.25pt}}%
{\ooalign{%
\ensuremath{%
\begin{picture}(8,8)
\thinlines
\put(1.5,2.75){\line(1,0){5}}
\end{picture}}%
\crcr
\hss\ensuremath{\textstyle\int}\hss}
\hspace*{-1.50pt}}%
{\ooalign{%
\ensuremath{%
\begin{picture}(8,8)
\thinlines
\put(2.5,2.0){\line(1,0){3}}
\end{picture}}%
\crcr
\hss\ensuremath{\scriptstyle\int}\hss}
\hspace*{-3.20pt}}%
{\ooalign{%
\ensuremath{%
\begin{picture}(8,8)
\thicklines
\put(2.98,1.5){\line(1,0){2}}
\end{picture}}%
\crcr
\hss\ensuremath{\scriptscriptstyle\int}\hss}
\hspace*{-3.0pt}}%
}

\numberwithin{equation}{section}

\def\notin{\not\in}

\begin{document}

\title[Local cut points and Ricci curvature]
{local cut points and metric measure spaces with Ricci curvature bounded below}
\author[Masayoshi Watanabe]{MASAYOSHI WATANABE}
\address{Mathematical Institute, Tohoku University, Sendai 980-8578, JAPAN}
\email{sa3m33@math.tohoku.ac.jp}
\thanks{This work was partially supported by
Research Fellowships of the Japan Society for the Promotion of Science for Young Scientists.}
\dedicatory{}
\date{\today}

\maketitle


\setlength{\baselineskip}{5mm}

\begin{abstract}
A local cut point is by definition 
a point that disconnects its sufficiently small neighborhood.
We show that
there exists an upper bound for the degree of a local cut point in a metric measure space
satisfying the generalized Bishop--Gromov inequality.
As a corollary,
we obtain an upper bound for the number of ends of such a space.
We also obtain some obstruction conditions for the existence of a local cut point
in a metric measure space satisfying the Bishop--Gromov inequality
or the Poincar\'{e} inequality.
For example,
the measured Gromov--Hausdorff limits of 
Riemannian manifolds with a lower Ricci curvature bound
satisfy these two inequalities.
\end{abstract}

\section{Introduction}

A point $x$ in a metric space is called a \textit{local cut point} if 
$U\setminus \{ x\}$ is disconnected for some connected neighborhood $U$ of $x$.
There exists no local cut point in 
any $n$-dimensional Alexandrov space
with curvature bounded from below for $n\ge 2$.
Every Gromov--Hausdorff limit 
of a sequence of Riemannian manifolds
with a uniform lower bound on sectional curvature
is such an Alexandrov space.
It is therefore natural to ask whether 
there exists a local cut point in the limit of manifolds 
with a uniform lower bound on Ricci curvature
unless the limit is one-dimensional.
We conjecture that the limit has no local cut point.

In this paper,
we consider metric measure spaces with ``Ricci curvature bounded below''. 
Let $\deg(x)$ denote the \textit{degree} of a point $x$, or
the supremum of the number of connected components of $U\setminus \{ x\}$
for all connected neighborhoods $U$ of $x$.
We give an upper bound for the degree of a local cut point.
As a consequence,
we obtain an upper bound for the number of ends
(see Subsection~\ref{ends} for the definition of an end).
We also obtain some obstruction conditions for the existence of a local cut point.

Cheeger and Colding~\cite{CCwarped}, \cite{CCI}, \cite{CCII}, \cite{CCIII},
and Menguy~\cite{M2}, \cite{M3}, \cite{M4} studied the limits of manifolds
with Ricci curvature bounded from below.
They constructed examples of limit spaces showing that
the local structure 
is more complicated than that of Alexandrov spaces.
Recently, Lott and Villani~\cite{LVricci}, \cite{LVweak}, 
Sturm~\cite{S}, \cite{SII}, and Ohta~\cite{O} independently introduced
the concept of lower bounds on Ricci curvature for metric measure spaces:
$N$-Ricci curvature $\ge K$, 
the curvature-dimension condition CD$(K,N)$, 
and the measure contraction property MCP$(K,N)$, respectively.
For each $N$-dimensional Riemannian manifold $M$, 
these three concepts are equivalent to that the Ricci curvature of $M$ is bounded below by $K$.
Moreover, these are preserved under the measured Gromov--Hausdorff limits.

We note that
the measured Gromov--Hausdorff limits of manifolds with Ricci curvature bounded from below
satisfy the generalized Bishop--Gromov inequality
and a Poincar\'{e} inequality of type $(1,1)$, as explained below.\\

Let $(X,d,\mu)$ be a complete, locally compact length space
equipped with a Borel measure.
We denote by $B_r(x)$ the open ball of radius $r$ and centered at $x\in X$.
Let $k\in \mathbb{R}$, $n\in \mathbb{N}$, $C\ge 1$, $0<r\le R$, and $x\in X$.
If $C=1$, then the following (\ref{usualBG}) is the usual Bishop--Gromov inequality
with respect to lower bound $(n-1)k$ of Ricci curvature 
and upper bound $n$ of dimension:
\begin{align}\label{usualBG}
\frac{\mu (B_R(x))}{\mu (B_r(x))}
\le C\,\frac{V_{k,\,n}(R)}{V_{k,\,n}(r)},
\end{align}
where $V_{k,\,n}(r)$
is the volume of a ball of radius $r$ in the $n$-dimensional, complete,
simply connected space of constant curvature $k$.

We consider a stronger inequality,
which is the directionally restricted version of (\ref{usualBG})
(\cite[Section 5. I $\! \!_+$]{G}, \cite[Appendix 2, (A.2.2.)]{CCI});
we call the inequality
the \textit{generalized Bishop--Gromov inequality with constant} $C$
(BG$(k,n)$ with $C$ for short).
This inequality is naturally extended to the case that $n\in \mathbb{R}$ with $n\ge 1$.
See Definition~\ref{defBG} for the precise definition
and \cite{KS} for slight different definitions.
This inequality is preserved under the measured Gromov--Hausdorff limits
(\cite[Section 5. I $\! \!_+$]{G}, \cite[Theorem 1.6, Theorem 1.10]{CCI}).
The measure contraction property
MCP$((n-1)k,n)$
implies BG$(k,n)$ with $C=1$.
For example,
the following metric measure spaces satisfy MCP$((n-1)k,n)$.
\begin{list}{}{}
  \item[\textbullet] $n$-dimensional Riemannian manifolds with Ricci curvature bounded below
                     by $(n-1)k$ equipped with Riemannian measure
  \item[\textbullet] $n$-dimensional Alexandrov spaces with curvature bounded below by $k$
                     equipped with the $n$-dimensional Hausdorff measure 
  \item[\textbullet] ``Nonbranching'' metric measure spaces satisfying 
                     the curvature-dimension condition CD$((n-1)k,n)$\\
                     See Definition~\ref{nonbranching} for the definition of nonbranching.
\end{list}
The non-Euclidean, finite-dimensional, normed linear spaces equipped with the Lebesgue measure
satisfy MCP$(0,n)$.
These spaces, however, can not arise as a Gromov--Hausdorff limit 
of any Riemannian manifolds with a uniform lower bound on Ricci curvature 
(see Proposition~\ref{nonlimit}).

Applying the method of the proof of Theorem 5.1 in \cite{CCII},
we have the following:

\begin{thm}\label{A}
Let $(X,d,\mu)$ be a metric measure space satisfying 
the generalized Bishop--Gromov inequality \textup{BG}$(k,n)$ with constant $C$
for some $k\in \mathbb{R}$, $n\ge 1$, and $C\ge 1$.
Assume that there exists a local cut point $x$ in $X$.
Then we have $\deg(x)\le C^2+1$.
\end{thm}

The Cheeger--Gromoll splitting theorem \cite{CG} states that, 
if a Riemannian manifold $M$ of nonnegative Ricci curvature contains a line,
then $M$ is isometric to $\mathbb{R}\times N$ for some manifold $N$.
Cheeger and Colding \cite{CCwarped} extended this 
to limit spaces of nonnegative Ricci curvature
in a generalized sense; see Theorem~\ref{splitting} of this paper.
For the limits of manifolds with Ricci curvature bounded from below ($C=1$),
Theorem~\ref{A} ($\deg(x)=2$) is also proved by using 
the Cheeger--Colding splitting theorem;
see Proposition~\ref{splitting-deg} of this paper.

By the splitting theorem,
a space of nonnegative Ricci curvature
has at most two ends.
We note that the splitting theorem for metric measure spaces satisfying
\textup{BG}$(0,n)$ with $C$ does not necessarily hold.
As a corollary of Theorem~\ref{A},
we have the following:

\begin{cor}\label{end}
Let $(X,d,\mu)$ be a metric measure space satisfying 
\textup{BG}$(0,n)$ with $C$
for some $n\ge 1$ and $C\ge 1$.
Then the number of ends of $X$ is at most $C^2+1$. 
\end{cor}

We investigate the geometric structure of 
the neighborhood of a local cut point (Theorem~\ref{diam}).
As a result, we see that the ``three-pronged'' space (see Figure~\ref{notBG})
does not satisfy BG$(k,n)$ with $C$.
We also study the structure of the accumulation of local cut points; we obtain that
the convergent sequence of certain local cut points 
``stands in a line'' (Corollary~\ref{standline}).\\

We now assume that a measure $\mu$ satisfies
$0<\mu(B_r(x))<+\infty$ for all $x\in X$ and all $0<r<+\infty$.
Let $1\le p<\infty$.
We say that a metric measure space $(X,d,\mu)$ satisfies
a \textit{Poincar\'{e} inequality of type $(1,p)$}, if
for all $R>0$ there exists a constant $C_P=C_P(p,R)>0$ depending only on $p$ and $R$ such that
\begin{equation}\label{introPoincare}
\mint_{B_r(x)} | u-u_{B_r(x)}|\ d\mu 
\le C_P\,r
\bigg(\mint_{B_r(x)} g^p\ d\mu\bigg)^{1/p}
\end{equation}
holds for all $x\in X$, $0<r\le R$, 
all measurable functions $u$, and all ``upper gradients'' $g$ of $u$ (see Section 4),
where 
$u_B:=\mint_Bu\ d\mu:=\mu(B)^{-1}\int_Bu\ d\mu$.

It follows from H\"{o}lder's inequality that
each metric measure space satisfying a Poincar\'{e} inequality of type $(1,p)$
also satisfies that of type $(1,q)$ 
for all $q\ge p$.
For instance, the following metric measure spaces satisfy 
a Poincar\'{e} inequality of type $(1,1)$.

\begin{list}{}{}
  \item[\textbullet] $n$-dimensional Riemannian manifolds 
                     with Ricci curvature bounded below by $(n-1)k$
                     equipped with Riemannian measure
                     (Buser \cite{B})\\
                     The constant $C_P$ in (\ref{introPoincare}) depends on $k,n$, and $R$.
                     If $k=0$, then $C_P$ depends only on $n$. See \cite[Theorem 5.6.5]{S-C}.
  \item[\textbullet] The measured Gromov--Hausdorff limits of 
                     manifolds with Ricci curvature bounded from below
                     (Cheeger and Colding {\cite[Theorem 2.15]{CCIII}})           
  \item[\textbullet] $n$-dimensional Alexandrov spaces equipped with 
                     the $n$-dimensional Hausdorff measure
                     (Kuwae, Machigashira, and Shioya {\cite[Theorem 7.2]{KMS}})
  \item[\textbullet] Nonbranching metric measure spaces satisfying
                     CD$(K,N)$ (von Renesse~\cite[Corollary in page 8]{R},
                     combined with Sturm~\cite[Lemma 4.1]{SII}) 
  \item[\textbullet] Metric measure spaces with a doubling measure and
                     the ``democratic condition'' 
                     DM (Lott and Villani~\cite[Theorem 2.5]{LVweak})                 
\end{list}
The space 
$\{ (x_1,x_2,\dots,x_n)\in \mathbb{R}^n\, |\, x_1^2+x_2^2+\cdots +x_{n-1}^2\le x_n^2\}$
equipped with the Euclidean distance and the $n$-dimensional Lebesgue measure
satisfies a Poincar\'{e} inequality of type $(1,p)$ for all $p>n$
(\cite[Example 4.2]{HajK}).

We obtain an obstruction condition for the existence of a local cut point 
as follows:

\begin{thm}\label{B}
Let $(X,d,\mu)$ be a metric measure space satisfying
a Poincar\'{e} inequality of type $(1,p)$ for some $1\le p<\infty$. 
Assume that
\begin{align}\label{katei}
\liminf_{r\to 0} \frac{\mu(B_r(x))}{r^p}=0
\end{align}
for a point $x\in X$.
Then $x$ is not a local cut point.
\end{thm}

Although the measured Gromov--Hausdorff limits of manifolds 
with Ricci curvature bounded from below
satisfy a Poincar\'{e} inequality of type $(1,1)$,
they do not necessarily satisfy 
the assumption (\ref{katei}) for $p=1$
(see \cite[Proposition 1.22]{CCI}).
Theorem~\ref{B} is one of the geometric consequences of the Poincar\'{e} inequality.
See \cite[Proposition 4.5]{HajK} and \cite[Theorem 3.13, Theorem 6.15]{HeiK} for related results.\\

The organization of this paper is as follows:
In Section 2, we recall the definitions of
the Hausdorff dimension, length spaces, and the (pointed, measured) 
Gromov--Hausdorff convergence.
In Section 3 we precisely define BG$(k,n)$ with $C$ and give its basic properties.
We prove Theorem~\ref{A} and Corollary~\ref{end} in Section 4.
Furthermore, we investigate the geometric structure of the neighborhood of local cut points,
and also study the structure of the accumulation of local cut points.
We prove Theorem~\ref{B} in Section 5.

\begin{ack}\upshape
The author would like to thank Prof.~Takashi Shioya
for his constant encouragement and valuable discussions.
He thanks Prof.~Tobias Holck Colding for useful comments.
He is indebted to Prof.~Koichi Nagano for stimulating his interest
in geometry of metric spaces.
He also thanks Prof.~Shin-ichi Ohta and Mr.~Shouhei Honda 
for valuable suggestions related to Theorem~\ref{B}.
\end{ack}

\section{Preliminaries}

In this section we first recall the definition of the Hausdorff dimension.
We then define length spaces and
the (pointed, measured) Gromov--Hausdorff convergence.
Let $(X,d)$ be a metric space.

\subsection{Hausdorff dimension}\label{Hausdorffdim}

We refer to \cite[Chapter 2]{AT} for details.
Let $A\subset X$.
For $0\le s<\infty$ and $0<\delta \le \infty$, we define
$$\mathcal{H}_{\delta}^s(A)=\omega_s\,
\inf \biggl\{ \sum_{i}\biggl( \frac{\mathrm{diam}(U_i)}{2}\biggr)^s
\biggm| A\subset \bigcup_{i}U_i,\, \mathrm{diam}(U_i)\le \delta \biggr\},$$
where 
$\omega_s:=\pi^{s/2}/\Gamma(s/2+1)$
and $\Gamma$ is the gamma function
$\Gamma(s)=\int_0^{\infty}e^{-x}x^{s-1}dx$.
If $s$ is a positive integer, 
then $\omega_s$ equals the volume of a unit ball in $\mathbb{R}^s$.
Note that $\mathcal{H}_{\delta}^s(A)$ is nonincreasing in $\delta$.
We define the $s$-\textit{dimensional Hausdorff measure} of $A$ by
$$\mathcal{H}^s(A)=\lim_{\delta \to 0}\mathcal{H}_{\delta}^s(A)
=\sup_{\delta>0}\mathcal{H}_{\delta}^s(A).$$
Then $\mathcal{H}^s$ is a Borel regular measure on $X$.
The \textit{Hausdorff dimension} of $A$ is defined as
$$\dim_{H}(A)=\inf \{ 0\le s<\infty \, |\, \mathcal{H}^s(A)=0\}
=\sup \{ 0\le s<\infty \, |\, \mathcal{H}^s(A)=+\infty\}.$$

\subsection{Length spaces}

For a continuous path $\gamma:[0,l]\to X$,
we define the length of $\gamma$ by
$$L(\gamma)=\sup_{0=t_0<t_1<\cdots <t_k=l}\sum_{i=1}^{k}
d(\gamma(t_{i-1}),\gamma(t_i)),$$
where the supremum is taken over all partitions of $[0,l]$.
By the triangle inequality, we have $L(\gamma)\ge d(\gamma(0),\gamma(l))$.
We say that a path $\gamma :[0, l]\to X$ is  
a \textit{geodesic} if it is locally minimizing and is proportional to arclength,
that is,  
for each $s\in [0,l]$ there exists $\epsilon = \epsilon(s)>0$ such that
\begin{align}\label{geod}
d(\gamma (t), \gamma (t^{\prime})) = |t - t^{\prime}|
\end{align}
holds for all $t, t^{\prime}\in (s-\epsilon , s+\epsilon)\cap [0, l]$.
Furthermore, we say that a path $\gamma:(-\infty,\infty)\to X$ is a \textit{line} if
(\ref{geod}) holds for all $t,t^{\prime}\in (-\infty,\infty)$.
In this paper, we assume that every path is proportional to arclength.

We say that $X$ is a \textit{length space} if 
$d(x,y)=\inf_{\gamma}L(\gamma)$ for all $x,y\in X$,
where the infimum is taken over all paths joining $x$ and $y$.
A metric space $X$ is a length space 
if and only if
for all $x,y\in X$ and all $\epsilon>0$ there exists a point $z\in X$ such that
$\max \{ d(x,z),d(z,y)\}\le d(x,y)/2+\epsilon$.
If $X$ is a complete, locally compact length space,
then all two points in $X$ are joined by a minimal geodesic.

See \cite{BBI}, \cite{G} and references therein for further information.

\subsection{(Pointed) Gromov--Hausdorff convergence}\label{GH}

Let us recall that the \textit{Hausdorff distance} between 
two closed bounded subsets $A$ and $B$ in a metric space $X$ is defined by
$$d_{H}(A,B)=\inf\{ \epsilon>0\, | \,A\subset U_{\epsilon}(B),B\subset U_{\epsilon}(A) \},$$
where $U_{\epsilon}(A)$ is the $\epsilon$-neighborhood of $A$.
Let $\mathcal{C}$ denote the set of isometry classes
of all compact metric spaces.
For $X,Y\in \mathcal{C}$,
the \textit{Gromov--Hausdorff distance} between $X$ and $Y$ is defined by
\begin{align*}
d_{GH}(X,Y)=
\inf \big\{ d_{H}(f(X),g(Y)) \bigm|\, &\text{all metric spaces}\ Z\ \mbox{and}\\
&\mbox{all isometric embeddings}\ f:X\to Z,\ g:Y\to Z\big\}.
\end{align*}
Then $d_{GH}$ defines a metric on $\mathcal{C}$.
Let $(X,d_X)$, $(Y,d_Y)$ be metric spaces. 
We say that for $\epsilon>0$
a map $\varphi:X\to Y$ is an $\epsilon$-\textit{approximation} if
the following two conditions hold:
\begin{itemize}
  \item[(i)]  $|d_X(x,y)-d_Y(\varphi(x),\varphi(y))|
                < \epsilon$ holds for all $x,y\in X$;
  \item[(ii)] the $\epsilon$-neighborhood of $\varphi(X)$ coincides with $Y$.
\end{itemize}
Let $X_i,X\in \mathcal{C}$ ($i=1,2,\dots$).
The sequence $\{ X_i\}$ Gromov--Hausdorff converges to $X$
($d_{GH}(X_i,X)\to 0$) as $i\to \infty$ if and only if
there exist $\epsilon_i$-approximations from $X_i$ to $X$
(or from $X$ to $X_i$) for some $\epsilon_i\to 0$.
Let $(X_i,x_i)$, $(X,x)$, $i=1,2,\dots$, be pointed metric spaces.
We say that $\{ (X_i,x_i)\}$ \textit{pointed Gromov--Hausdorff converges} to $(X,x)$, if
for each $R>0$ 
there exist $R_i\searrow R$, $\epsilon_i\searrow 0$, 
and $\epsilon_i$-approximations $\varphi_i:B_{R_i}(x_i)\to B_R(x)$
with $\varphi_i(x_i)=x$.
A pointed Gromov--Hausdorff limit of $\{ (X,r_i^{-1}d,x)\}$ as $r_i\to 0$
is called a \textit{tangent cone} at $x$.
The Gromov--Hausdorff limit of a sequence of length spaces is also a length space.

We refer to \cite{BBI}, \cite{G} for details.\\

Let $(M_i,p_i)$, $i=1,2,\dots$, be $n$-dimensional, complete pointed Riemannian manifolds
with Ricci curvature $\mathrm{Ric}_{M_i}\ge -(n-1)$.
Gromov's compactness theorem states that
$\{(M_i,p_i)\}$ pointed Gromov--Hausdorff subconverges to some pointed metric space $(X,x)$.
We then have $\dim_{H}(X)\le n$ (see Proposition~\ref{dim}).
The sequence $\{ M_i\}$ is said to \textit{collapse} to $X$ if $\dim_{H}(X)<n$.
In fact, if $\{ M_i\}$ collapses, then $\dim_{H}(X)\le n-1$ (\cite[Theorem 3.1]{CCI}).
It is conjectured that the Hausdorff dimension of the limit $X$ 
is an integer (Fukaya's conjecture~\cite[Conjecture 3.13]{F}).  
By Gromov's compactness theorem,
for every point $x$ in the limit $X$
there exists a tangent cone at $x$. 
A tangent cone at each point is not necessarily unique
even in the noncollapsed case
(\cite[Example 8.41]{CCI}).
It does not always have a metric cone structure in the collapsed case
(\cite[Example 8.95]{CCI}).
We state the splitting theorem for the limit,
which holds even in the collapsed case.

\begin{thm}[Cheeger--Colding~\cite{CCwarped}]\label{splitting}
Let $(X,x)$ be a pointed Gromov--Hausdorff limit of complete Riemannian manifolds $M_i$
with $\mathrm{Ric}_{M_i}\ge -\delta_i$, where $\delta_i\to 0$.
If $X$ contains a line,
then $X$ is isometric to $\mathbb{R}\times Y$ for some length space $Y$.
\end{thm}

\subsection{Measured Gromov--Hausdorff convergence}

Fukaya~\cite{F} introduced
the concept of the measured Gromov--Hausdorff convergence.
Let $\mathcal{CM}$ denote the set of all compact metric spaces $X$
equipped with a Borel measure $\mu$ such that $\mu(X)\le 1$.
A directed system 
$\{ (X_{\alpha},\mu_{\alpha})\}_{\alpha \in \mathcal{A}}\subset \mathcal{CM}$
is said to \textit{measured Gromov--Hausdorff converge} to $(X,\mu)\in \mathcal{CM}$ if
for each $\alpha \in \mathcal{A}$
there exist $\epsilon_{\alpha}>0$ and
a Borel measurable $\epsilon_{\alpha}$-approximation 
$\varphi_{\alpha}:X_{\alpha}\to X$ such that
\begin{itemize}
  \item[(i)]  $\lim_{\alpha \in \mathcal{A}}\epsilon_{\alpha}=0;$
  \item[(ii)] a directed system of push-forward measures
              $\{ (\varphi_{\alpha})_*\mu_{\alpha}\}_{\alpha}$
              converges to $\mu$ vaguely:
              $$\lim_{\alpha \in \mathcal{A}}
              \int_{X_{\alpha}}f\circ \varphi_{\alpha} \ d\mu_{\alpha}
              =\int_Xf\ d\mu$$
              holds for all continuous functions $f:X\to \mathbb{R}$.
\end{itemize}
We induce a topology on $\mathcal{CM}$
by the measured Gromov--Hausdorff convergence. 
The topology is Hausdorff (\cite[Proposition 2.7]{F}),
and the projection from $\mathcal{CM}$ to $\mathcal{C}$ is proper (\cite[Proposition 2.10]{F}).
We refer to \cite{G} for other topologies on $\mathcal{CM}$.\\

Let $M_i$, $i=1,2,\dots$, be compact $n$-dimensional Riemannian manifolds
with $\mathrm{Ric}_{M_i}\ge -(n-1)$.
Assume that $\{ M_i\}$ Gromov--Hausdorff converges to
a compact metric space $(X,d)$.
Since $(M_i,\mathrm{dvol}_{M_i}/\mathrm{Vol}(M_i))\in \mathcal{CM}$,
where $\mathrm{dvol}_{M_i}$ is the volume element of $M_i$,
there exists a measure $\mu$ on $X$ such that 
$\{ (M_i,\mathrm{dvol}_{M_i}/\mathrm{Vol}(M_i))\}$
measured Gromov--Hausdorff subconverges to $(X,\mu)$.
In the noncollapsed case,
$\mu$ is unique and it coincides with the $n$-dimensional Hausdorff measure,
up to constant multiple (\cite[Theorem 5.9]{CCI}).
In the collapsed case,
$\mu$ is not necessarily unique (\cite[Example 1.24, Section 8]{CCI}).

\section{The generalized Bishop--Gromov inequality}

In this section,
we define the generalized Bishop--Gromov inequality 
BG$(k,n)$ with constant $C$,
and give basic properties of 
a metric measure space satisfying BG$(k,n)$ with $C$. 

We denote by $c(c_1,\dots,c_l)$
a positive constant depending only on $c_1,\dots,c_l$.
Let $(X,d,\mu)$ be a complete, locally compact length space
equipped with a Borel measure.

\subsection{The definition of the generalized Bishop--Gromov inequality}

Let $k, n\in \mathbb{R}$ with $n\ge 1$.
For $0\le r_1<r_2$,
we define
$$V_{k,\,n}(r_1,r_2)=\alpha_{n-1}\int_{r_1}^{r_2}\mathsf{s}_{k}(t)^{n-1}\ dt,$$
where 
$\alpha_{n-1}:=2\pi^{n/2}/\Gamma(n/2)=n\,\omega_n$ and
$\mathsf{s}_{k}(t)$, $t\ge 0$, is defined by 
\begin{equation*}
  \mathsf{s}_{k}(t)=\begin{cases} 
           \frac{1}{\sqrt{k}}\sin ({\sqrt{k}\,t})     & \mbox{if}\ k>0,\\                  
           t                                          & \mbox{if}\ k=0,\\
           \frac{1}{\sqrt{-k}}\sinh ({\sqrt{-k}\,t})  & \mbox{if}\ k<0.
         \end{cases}
\end{equation*}
If $n$ is a positive integer, then
$V_{k,\,n}(r_1,r_2)$ is equal to the volume of an annulus
of radius between $r_1$ and $r_2$ in the $n$-dimensional, complete,
simply connected space of constant curvature $k$.
We denote by $B_r(x)$ and $\overline{B}_r(x)$ the open and closed ball
of radius $r$ and centered at $x$, respectively.
For a point $x\in X$,
let $U$ be a measurable set in $A_{r_1,\,r_2}(x)$, where  
$A_{r_1,\,r_2}(x)=B_{r_2}(x)\setminus \overline{B}_{r_1}(x)$
and $A_{0,\,r}(x)=B_{r}(x)$.
For $0\le s_1<s_2$ with $s_1\le r_1$ and $s_2\le r_2$, set
$$S_{s_1,\,s_2}(x,U)=\big\{ y\in A_{s_1,\,s_2}(x) \bigm| 
                     d(x,y)+d(y,z)=d(x,z) \mbox{ \textup{for some} } z\in U\big\}$$
(that is, 
$S_{s_1,\,s_2}(x,U)$ is the intersection of $A_{s_1,\,s_2}(x)$ and
all minimal geodesics between $x$ and each point in $U$).

\begin{defn}\label{defBG}\upshape
Let $C\ge 1$. 
We say that a metric measure space $(X,d,\mu)$
satisfies the \textit{generalized Bishop--Gromov inequality with constant} $C$
(BG$(k,n)$ with $C$ for short), if
for all points $x\in X$
and all $0\le r_1<r_2$, $0\le s_1<s_2$ with $s_1\le r_1, s_2\le r_2$,
\begin{align}\label{BG(k,n)}
\frac{\mu (U)}{\mu (S_{s_1,\,s_2}(x,U))}
\le C\, \frac{V_{k,\,n}(r_1,r_2)}{V_{k,\,n}(s_1,s_2)}
\end{align}
holds for all measurable sets $U$ in $A_{r_1,\,r_2}(x)$.
\end{defn}

BG$(k,n)$ with $C=1$ induces the usual one.
Indeed, if we choose $r_1=s_1=0$, $s_2=r\le R=r_2$, and $U=B_R(x)$,
then (\ref{BG(k,n)}) implies
$$\frac{\mu (B_R(x))}{\mu (B_r(x))}
\le \frac{V_{k,\,n}(0,R)}{V_{k,\,n}(0,r)}.$$

\subsection{Examples}

We recall metric measure spaces satisfying BG$(k,n)$ with $C=1$,
as shown in the introduction.
Every $n$-dimensional, normed linear space
equipped with the $n$-dimensional Lebesgue measure
satisfies BG$(0,n)$ with $C=1$.

\begin{prop}\label{nonlimit}
Let $X$ be an $n$-dimensional, normed linear space.
Assume that $X$ is a Gromov--Hausdorff limit of Riemannian manifolds $M_i$ 
with Ricci curvature bounded from below.
Then $X$ is isometric to $\mathbb{R}^n$.
\end{prop}

\begin{proof}
Expanding the metrics of converging manifolds, we may assume that 
$\mathrm{Ric}_{M_i}\ge -\delta_i$, where $\delta_i\to 0$.
Since $X$ contains $n$ orthogonal lines,
the splitting theorem (Theorem~\ref{splitting}) completes the proof.
\end{proof}

\begin{defn}[{\cite[Definition 2.8]{S}}]\label{nonbranching}\upshape
A metric space $(X,d)$ is said to be \textit{nonbranching} if
for all four points $y,x_0,x_1,x_2\in X$
such that $y$ is a midpoint between $x_0$ and $x_1$ and between $x_0$ and $x_2$,
we have $x_1=x_2$.
A point $y$ is called a \textit{midpoint} between $x_0$ and $x_1$
if $2^{-1}d(x_0,x_1)=d(x_0,y)$ holds.
\end{defn}
A metric space $X$ is nonbranching if and only if
for all two minimal geodesics $\gamma,\gamma^{\prime}:[0,l]\to X$ 
with $\gamma(0)=\gamma^{\prime}(0)$,
we have $\inf \{ t>0\,|\,\gamma(t)\neq \gamma^{\prime}(t)\}=0$ or $\infty$,
where $\inf \emptyset:=\infty$.
We see that each Alexandrov space with curvature bounded from below
is nonbranching.

\begin{rem}\upshape
For nonbranching metric measure spaces, 
the curvature-dimension condition CD$(K,N)$ implies 
the measure contraction property MCP$(K,N)$
(see \cite[Lemma 4.1, Theorem 5.4]{SII}).
Moreover, CD$(K,N)$ is a little weaker condition than having 
$N$-Ricci curvature $\ge K$.

All $n$-dimensional Alexandrov spaces with curvature bounded below by $k$
equipped with the $n$-dimensional Hausdorff measure $\mathcal{H}^n$
satisfy MCP$((n-1)k,n)$ 
(\cite[Proposition 2.8]{O}, \cite[Theorem 5.7]{SII}, \cite[Lemma 6.1]{KS}).
For Riemannian manifolds,
lower bounds of sectional curvature
imply that of Ricci curvature.
However, it is an open problem whether those Alexandrov spaces satisfy CD$((n-1)k,n)$ 
or have $n$-Ricci curvature $\ge (n-1)k$.
In particular,
we do not know whether all Alexandrov spaces
can arise as Gromov--Hausdorff limits of Riemannian manifolds 
with Ricci curvature bounded from below.
\end{rem}

\begin{rem}\upshape
Let $(X,d,\mu)$ be a measured Gromov--Hausdorff limit of 
$n$-dimensional compact Riemannian manifolds $M_i$ with $\mathrm{Ric}_{M_i}\ge (n-1)k$
and equipped with the normalized Riemannian measure. 
Set $m=\dim_{H}(X)$, which is not necessarily an integer.
It is conjectured that there exists a number $l$ with $m\le l\le n$ such that
$(\ref{BG(k,n)})$ holds for
$\mu=\mathcal{H}^m$, $n=l$, and $C=1$ (\cite[Conjecture 1.34]{CCI}).
\end{rem}

\subsection{Basic properties}\label{basic}

We begin by recalling the definition of a doubling measure.
We say that a measure $\mu$ is \textit{doubling} if
for all $R>0$ 
there exists a constant $C_D=C_D(R)\ge1$ such that $\mu(B_{2r}(x))\le C_D\, \mu(B_r(x))$
holds for all $x\in X$ and all $0<r\le R$.
If $\mu$ is doubling, 
then $X$ is proper $($that is, 
all closed bounded subsets are compact$)$.
We note that if $(X,d,\mu)$ satisfies BG$(k,n)$ with $C$,
then $\mu$ is doubling.

Although the following lemmas and proposition are somewhat standard,
we prove them for the completeness of this paper.

\begin{lem}
Let $(X,d,\mu)$ satisfy 
\textup{BG}$(k,n)$ with $C=1$ for some $k,n\in \mathbb{R}$ with $n\ge 1$.
Then we have
$\mu(S_r(x))=0$ for all $x\in X$ and all $r>0$,
where $S_r(x)$ is the sphere of radius $r$ and centered at $x$.
In particular, 
if $X$ does not consist of a single point,
then $\mu(\{ x\})=0$ for all $x\in X$.
\end{lem}

\begin{proof}
For a sufficiently small $\epsilon>0$,
we have
\begin{align*}
\mu(S_r(x))
&\le \mu(B_{r+\epsilon}(x))-\mu(B_{r-\epsilon}(x))
= \mu(B_{r-\epsilon}(x))
\bigg( \frac{\mu(B_{r+\epsilon}(x))}{\mu(B_{r-\epsilon}(x))}-1\bigg)\\
&\le \mu(B_{r-\epsilon}(x))
\bigg( \frac{V_{k,n}(0,r+\epsilon)}{V_{k,n}(0,r-\epsilon)}-1\bigg),
\end{align*}
where we have used (\ref{BG(k,n)}) for the last inequality.
Letting $\epsilon \to 0$ completes the proof. 
\end{proof}

Let $A\subset X$.
We say that $A$ is \textit{convex} if
for all two points $x,y\in A$,
all minimal geodesics between $x$ and $y$ are contained in $A$.

\begin{lem}\label{basiclem}
Let $(X,d,\mu)$ satisfy
\textup{BG}$(k,n)$ with $C$ for some $k\in \mathbb{R}$, $n\ge 1$, and $C\ge 1$.
Then the following holds:
\begin{itemize}
\item[(1)] If $k^{\prime}\le k$ and if $n^{\prime}\ge n$,
      then $(X,d,\mu)$ satisfies \textup{BG}$(k^{\prime},n^{\prime})$ with $C$.
\item[(2)] If $a,b>0$,
      then $(X,ad,b\mu)$ satisfies \textup{BG}$(k/a^2,n)$ with $C$.
\item[(3)] If $A$ is a convex subset in $X$,
      then $(A,d,\mu)$ also satisfies \textup{BG}$(k,n)$ with $C$.
\end{itemize}      
\end{lem}

\begin{proof}
Fix $0\le r_1<r_2,0\le s_1<s_2$ with $s_1\le r_1, s_2\le r_2$.
(1) follows from  facts that $V_{k,\,n}(r_1,r_2)/V_{k,\,n}(s_1,s_2)$ is 
monotone nonincreasing in $k$
and monotone nondecreasing in $n$.
We denote by $A_{r_1,\,r_2}^{(d)}(x)$ and $S_{s_1,\,s_2}^{(d)}(x,U)$ 
the sets $A_{r_1,\,r_2}(x)$ and $S_{s_1,\,s_2}(x,U)$ 
with respect to a metric $d$, respectively.
For all measurable sets $U$ in $A_{r_1,\,r_2}^{(ad)}(x)$, we have
$$\frac{(b\mu)(U)}{(b\mu)\big( S_{s_1,\,s_2}^{(ad)}(x,U)\big)}
=\frac{\mu(U)}{\mu\big( S_{s_1/a,\,s_2/a}^{(d)}(x,U)\big)}
\le C\,\frac{V_{k,\,n}(r_1/a,r_2/a)}
{V_{k,\,n}(s_1/a,s_2/a)}
=C\, \frac{V_{k/a^2,\,n}(r_1,r_2)}{V_{k/a^2,\,n}(s_1,s_2)},$$
which implies $(2)$.
Let $x\in A$ and $U\subset A_{r_1,r_2}(x)$.
Since $A$ is convex, $S_{s_1,\,s_2}(x,U)$ is contained in $A$.
We hence have $(3)$.
\end{proof}

Let $\epsilon>0$.
A set $S\subset X$ is called an $\epsilon$-\textit{net}
if we have $d(x,S)\le \epsilon$ for all $x\in X$.
$S$ is said to be $\epsilon$-\textit{separated}
if we have $d(x,y)\ge \epsilon$ for all two distinct points $x,y\in S$.
We see that each maximal $\epsilon$-separated set is an $\epsilon$-net.

\begin{prop}\label{dim}
Let $(X,d,\mu)$ be as in Lemma~\ref{basiclem}.
Then we have $\dim_{H}(X)\le n$.
\end{prop}

\begin{proof}
We will show that $\mathcal{H}^{n^{\prime}}(X)=0$ holds
for all $n^{\prime}>n$.
We may assume that $\mathrm{diam}(X)<\infty$.
Let $\{ x_i\}_{i=1}^N$ be a maximal $\delta$-separated set in $X$ for $\delta>0$.
Then $\{ x_i\}_{i=1}^N$ is a $\delta$-net.
Note that there exists a constant $c(k,n)>0$ depending only on $k$ and $n$
such that $V_{k,\,n}(0,\delta)\ge c(k,n)\delta^n$ holds
for all (sufficiently small if $k>0$) $\delta>0$.
Let $x_{i_0}$ be a point such that
$\mu(B_{\delta/2}(x_{i_0}))$ attains 
$\min_{1\le i\le N} \mu(B_{\delta/2}(x_i))$.
Remark that $\{ B_{\delta/2(x_i)}\}_{i=1}^N$ is disjoint.
Since $(X,d,\mu)$ satisfies BG$(k,n)$ with $C$,
we have
\begin{align*}
\mu(X)
\ge \sum_{i=1}^N\mu(B_{\delta/2}(x_i))
\ge N\mu(B_{\delta/2}(x_{i_0}))
\ge NC^{-1}\frac{V_{k,\,n}(0,\delta/2)}{V_{k,\,n}(0,\mathrm{diam}(X))}\,
\mu\big(B_{\mathrm{diam}(X)}(x_{i_0})\big).
\end{align*}
Hence, $N\le c(k,n)\delta^{-n}$, which implies
\begin{align*}
\min \bigg\{ N\in \mathbb{N}\biggm| X= 
\bigcup_{i=1}^NB_{\delta}(x_i) \bigg\}\\
&\hspace{-11em}\le \max \big\{ N\in \mathbb{N}\bigm| 
\{ x_i\}_{i=1}^N\mbox{ is a maximal\ } \delta\mbox{-separated set}\big\}\\
&\hspace{-11em}\le c(k,n)\delta^{-n}.
\end{align*}
Therefore,
\begin{align*}
\mathcal{H}_{\delta}^{n^{\prime}}(X)
&=\omega_{n^{\prime}}\,\inf \Biggl\{ \sum_{i=1}^{\infty}
\biggl( \frac{\mathrm{diam}(U_i)}{2}\biggr)^{n^{\prime}}
\Biggm| X= \bigcup_{i=1}^{\infty}U_i,\mathrm{diam}(U_i)\le \delta \Biggr\}\\
&\le \omega_{n^{\prime}}\,\inf \Biggl\{ \sum_{i=1}^{\infty}
\biggl( \frac{\mathrm{diam}(B_{\delta/2}(x_i))}{2}\biggr)^{n^{\prime}}
\Biggm| X= \bigcup_{i=1}^{\infty}B_{\delta/2}(x_i) \Biggr\}\\
&\le \omega_{n^{\prime}}\,\Big(\frac{\delta}{2}\Big)^{n^{\prime}}
\min \bigg\{ N\in \mathbb{N}\biggm| X= 
\bigcup_{i=1}^NB_{\delta/2}(x_i) \bigg\}\\
&\le c(k,n,n^{\prime})\,\delta^{n^{\prime}-n}.
\end{align*}
Letting $\delta \to 0$ completes the proof.
\end{proof}

\section{Local cut points}

In this section, we first prove Theorem~\ref{A}.
In Subsection~\ref{ends},
we recall the definition of an end,
and then prove Corollary~\ref{end}.
Furthermore, we investigate the geometric structure of the neighborhood of a local cut point,
and also study the structure of the accumulation of local cut points
by using the idea of the proof of Theorem~\ref{A}.
We now assume that $(X,d)$ is a complete, locally compact length space.

\subsection{Local cut points}

\begin{defn}[{local cut point~\cite[3.32]{G}, $r$-cut point}]\label{cutpt}\upshape
We say that a point $x\in X$ is a \textit{local cut point} if
$U\setminus \{ x\}$ is disconnected for some connected neighborhood $U$ of $x$.
The \textit{degree} of $x$, denoted by $\deg(x)$, is defined as
the supremum of the number of connected components of $U\setminus \{ x\}$
for all connected neighborhoods $U$ of $x$.
Let $r>0$.
We say that a point $x\in X$ is an $r$-\textit{cut point} if
the following three conditions hold:
\begin{itemize}
  \item[(i)]   $\overline{B}_r(x)\setminus \{ x \}$ is disconnected;
  \item[(ii)]  the number of connected components of $\overline{B}_r(x)\setminus \{ x \}$  
               is equal to $\deg(x)$;
  \item[(iii)] $O\cap S_r(x)$ is nonempty for all connected components $O$ of
               $\overline{B}_r(x)\setminus \{ x \}$.
\end{itemize}
\end{defn}

If $x$ is a local cut point,
then $U\setminus \{ x\}$ is disconnected 
for \textit{every} sufficiently small neighborhood $U$ of $x$.
The end points in a graph (one-dimensional space) 
are \textit{not} local cut points.
An interior point in a graph
is not always a local cut point; see Example~\ref{cutptex} (3).
We have $\deg(x)\ge 2$ for each local cut point $x$.
Every local cut point with finite degree is an $r$-cut point for some $r>0$.
Let $0<r_1\le r_2$.
If a point $x$ is an $r_2$-cut point, then $x$ is an $r_1$-cut point.

\begin{ex}\upshape \label{cutptex}
(1) Consider the set 
$\bigcup_{i=0}^{\infty}\{ r\exp(2^{-i}\pi \sqrt{-1}\,)\in \mathbb{C}\,|\,0\le r\le 2^{-i}\}$
with the induced distance (\cite[Example 9.1]{Sy}).
The origin is a local cut point and its degree is infinite.

(2) Consider the set
$\bigcup_{i=1}^{\infty}\{ (x,y)\in \mathbb{R}^2\,|\,x^2+(y-1/i)^2=1/i^2\}$
with the induced distance.
The origin is a local cut point and its degree is infinite.
 
(3) Consider the set
$\{ (x,0)\,|\,0\le x\le 1\}\cup \{ (0,y)\,|\, 0\le y\le 1\}
\cup \big( \bigcup_{i=0}^{\infty}\{ (x,-x+2^{-i})\,|\, 0\le x\le 2^{-i})\} \big)
\subset \mathbb{R}^2$
with the induced distance. 
The origin is \textit{not} a local cut point.

(4) Consider the interval $[0,1]\subset \mathbb{R}$ with the Euclidean distance.
The point $1/3\in [0,1]$ is \textit{not} an $r$-cut point for any $r>1/3$
because (iii) is not satisfied.

(5) Consider the set $\{ (x,y)\,|\,x^2+y^2=1\}\cup 
\{ (x,0)\,|\,1\le x\} \subset \mathbb{R}^2$
with the induced distance.
The point $(-1,0)$ is \textit{not} an $r$-cut point for any $r\ge \pi$
because (ii) is not satisfied.
\end{ex}

\begin{figure}[tbp]
\begin{center}
 \includegraphics[width=10cm,clip]{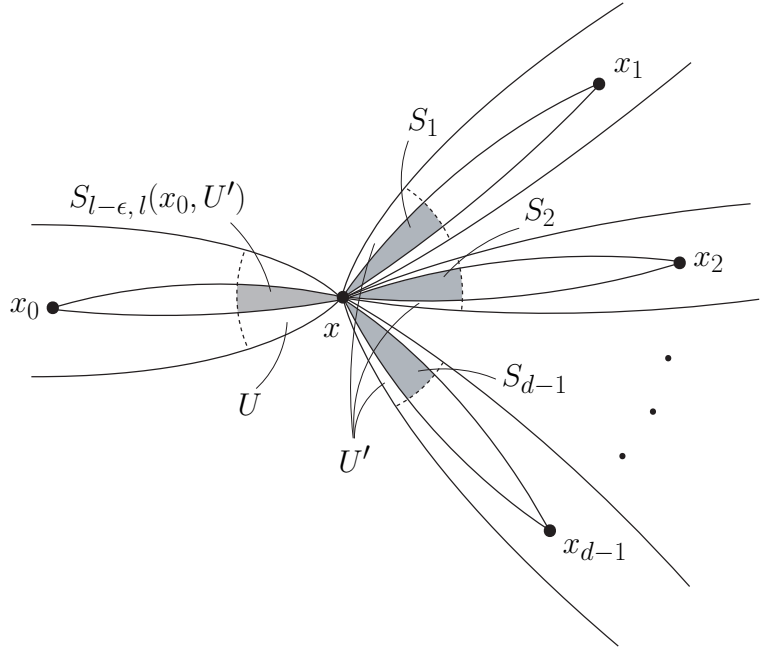}
\caption{Proof of Theorem~\ref{A}}
\label{ProofofThmA}
\end{center}
\end{figure}

\begin{proof}[Proof of Theorem~\ref{A}]
For any positive integer $d$ with $d\le \deg(x)$,
take a sufficiently small $r>0$.
Choose any $d$ connected components $O_1,O_2,\dots,O_d$
of $\overline{B}_r(x)\setminus \{ x\}$
such that 
the degree of $x$ in
$\overline{B}_r(x)\cap (O_1\cup \cdots \cup O_d)$
is equal to $d$, and
$O_i\cap S_r(x)$ is nonempty for all $1\le i\le d$.
Let $l$ be a positive number
with $0<l\le r/2$.
For each $1\le i\le d$,
we choose a point $x_i\in O_i$ 
such that $d(x,x_i)=l$.
See Figure~\ref{ProofofThmA}.

For a sufficiently small $0<\epsilon \ll l$,
set $U=B_{\epsilon}(x)\cap O_1$.
We claim that each minimal geodesics between every point in $U$ and $x_i$ ($2\le i \le d$)
passes through the local cut point $x$.
Suppose that there exists a minimal geodesic $\gamma:[0,l]\to X$
from some point $y$ in $U$ to $x_i$
such that $\gamma$ does not pass through $x$.
By the choice of $r$ and $l$,
the point $\gamma(t)$ is not contained in $B_r(x)$ for some $t\in [0,l]$.
Therefore,
$$d(y,x_i)=d(\gamma(0),\gamma(l))=l\ge t=d(y,\gamma(t))\ge d(x,\gamma(t))-d(x,y)>r-\epsilon.$$
On the other hand,
$$d(y,x_i)\le d(y,x)+d(x,x_i)<\epsilon+l\le \epsilon +\frac{r}{2}.$$
This is a contradiction.

We first use the generalized Bishop--Gromov inequality~(\ref{BG(k,n)})
with the base point $x_i$ ($2\le i\le d$).
We see that
$A_{l,\,l+\epsilon}(x_i)$ contains $U$
for $2\le i\le d$.
Set $S_i=S_{l-\epsilon,\,l}(x_i,U)$ for $2\le i\le d$.
We recall that $S_i$ is the intersection of $A_{l-\epsilon,\,l}(x_i)$
and all minimal geodesics from $x_i$ to each point in $U$, hence to $x$.
Then $S_i\cap S_j$ is empty for all $i\neq j$.
By applying (\ref{BG(k,n)})
to $s_1=l-\epsilon$, $s_2=r_1=l$, $r_2=l+\epsilon$, the point $x_i$, and the set $U$,
it follows for $2\le i\le d$ that
\begin{align}\label{degS_i}
C^{-1} \frac{V_{k,\,n}(l-\epsilon,l)}
     {V_{k,\,n}(l,l+\epsilon)}
\le \frac{\mu\big(S_{l-\epsilon,\,l}(x_i,U)\big)}{\mu (U)}
= \frac{\mu(S_i)}{\mu(U)}.
\end{align}
We now denote $U^{\prime}=\bigcup_{i=2}^d S_i$.
Since $\bigcap_{i=2}^{d}S_i=\emptyset$,
we have $\sum_{i=2}^{d}\mu(S_i)= \mu(U^{\prime})$.
Summing up (\ref{degS_i}) for all $2\le i\le d$, we obtain
\begin{align}\label{degS_1+S_2}
(d-1)\,C^{-1} \frac{V_{k,\,n}(l-\epsilon,l)}
     {V_{k,\,n}(l,l+\epsilon)}
\le \frac{\sum_{i=2}^{d}\mu(S_i)}{\mu (U)}
= \frac{\mu(U^{\prime})}{\mu(U)}.
\end{align}

Next, we use (\ref{BG(k,n)}) with the base point $x_1$.
We see that $A_{l,\,l+\epsilon}(x_1)$ contains $U^{\prime}$.
Applying (\ref{BG(k,n)})
to $s_1=l-\epsilon$, $s_2=r_1=l$, $r_2=l+\epsilon$, the point $x_1$, and the set $U^{\prime}$, 
we obtain
\begin{align}\label{degUprime}
\frac{\mu(U^{\prime})}
{\mu\big(S_{l-\epsilon,\,l}(x_1,U^{\prime})\big)}
\le 
C\, \frac{V_{k,\,n}(l,l+\epsilon)}
     {V_{k,\,n}(l-\epsilon,l)}.
\end{align}
Note that each minimal geodesics between $x_0$ and each point in $U^{\prime}$ 
passes through $x$
similarly to that mentioned above; hence,
$U$ contains $S_{l-\epsilon,\,l}(x_1,U^{\prime})$.
Combining (\ref{degS_1+S_2}) and (\ref{degUprime}), we have
\begin{align*}
d\le C^2\bigg( \frac{V_{k,\,n}(l,l+\epsilon)}{V_{k,\,n}(l-\epsilon,l)}\bigg)^2+1.
\end{align*}
Letting $\epsilon \to 0$, 
we have $d\le C^2+1$.
This completes the proof.
\end{proof}

For the limits of Riemannian manifolds with Ricci curvature bounded from below,
the Cheeger--Colding splitting Theorem (Theorem~\ref{splitting}) implies
the same conclusion in Theorem~\ref{A};
see Proposition~\ref{splitting-deg} below.
Since the splitting theorem for metric measure spaces 
satisfying BG$(0,n)$ with $C$ do not hold in general (see Proposition~\ref{nonlimit}),
the proof of Theorem~\ref{A} gives another proof that $\deg(x)=2$.

\begin{prop}\label{splitting-deg}
Let $(X,d)$ be a Gromov--Hausdorff limit of 
$n$-dimensional, complete Riemannian manifolds $(M_i,d_i)$
with $\mathrm{Ric}_{M_i}\ge -(n-1)$,
where $d_i$ is the Riemannian distance on $M_i$.
If there exists a local cut point $x$ in $X$,
then every tangent cone at $x$ is isometric to $\mathbb{R}$,
in particular $\deg(x)=2$.
\end{prop}

\begin{proof}
Take points $p_i$ in $M_i$ such that $\{(M_i,d_i,p_i)\}$ pointed Gromov--Hausdorff
converges to $(X,d,x)$.
Let us recall that a tangent cone at $x$ is 
the pointed limit space of $(X,r_i^{-1}d,x)$ as $r_i \to 0$.
By passing to a subsequence,
the tangent cone at $x$ is itself the pointed limit of 
rescaled manifolds $\big( M_j,r_{j}^{-1}d_j,p_j\big)$.
Then the Ricci curvature of $\big( M_j,r_{j}^{-1}d_j\big)$ is bounded below 
by $-(n-1)r_{j}^2$.

Because every local cut point is an interior point of some geodesic,
the tangent cone at $x$ contains a line (expanded from the geodesic).
Applying the splitting theorem (Theorem~\ref{splitting}),
we see that the tangent cone at $x$ is isometric to $\mathbb{R}\times Y$
for some length space $Y$. 
Since $x$ is a local cut point, it follows that $Y$ consists of a single point.
\end{proof}

\begin{rem}\upshape
We recall that 
if $(X,d,\mu)$ satisfies BG$(k,n)$ with $C$,
then $\mu$ is doubling.
Although the length space in Example~\ref{cutptex} (1) 
has the natural doubling measure,
the degree of the origin is infinite.
\end{rem}

In the case of graphs,
we obtain a better estimate of the degree.

\begin{prop}\label{graph}
Let $(X,d,\mathcal{H}^1)$ be a graph,
where $d$ is the usual distance 
and $\mathcal{H}^1$ is the one-dimensional Hausdorff measure.
Assume that $(X,d,\mathcal{H}^1)$ satisfies 
\textup{BG}$(k,n)$ with $C$ 
for some $k\in \mathbb{R}$, $n\ge 1$, and $C\ge 1$.
Then we have $\deg(x) \le C+1$ for all $x\in X$.
\end{prop}

\begin{proof}
Let $x$ be a local cut point (an interior point) in $X$.
Take a point $x_1\in X$ with $d(x,x_1)\ll 1$, 
and denote $l=d(x,x_1)$.
Let $\gamma:[0,l]\to X$ be the minimal geodesic from $x_1$ to $x$.
Set $U=B_{\epsilon}(x)\setminus \gamma([0,l])$ for $0<\epsilon \ll 1$.
We then have 
$S_{l-\epsilon,\,l}(x_1,U)=\gamma((l-\epsilon,l))$
as $l=d(x,x_1)$ is sufficiently small.
By using (\ref{BG(k,n)})
to $s_1=l-\epsilon$, $s_2=r_1=l$, $r_2=l+\epsilon$, the point $x_1$, and the set $U$,
it follows that
\begin{align}
\frac{(\deg(x)-1)\epsilon}{\epsilon}=
\frac{\mathcal{H}^1(U)}
{\mathcal{H}^1\big(S_{l-\epsilon,\,l}(x_1,U)\big)}
\le 
C\, \frac{V_{k,\,n}(l,l+\epsilon)}
     {V_{k,\,n}(l-\epsilon,l)}.
\end{align}
Taking $\epsilon \to 0$,
we obtain $\deg(x)\le C+1$.
\end{proof}

\begin{rem}\upshape
In general,
the converse of Proposition~\ref{graph} does not necessarily hold.
\end{rem}

\subsection{The number of ends}\label{ends}

We recall the definition of an end.

\begin{defn}\upshape
Let $\gamma_1,\gamma_2:[0,\infty)\to X$ be rays 
from the base point $x$.
Two rays $\gamma_1$ and $\gamma_2$ are said to be \textit{cofinal} if
$\gamma_1(t)$ and $\gamma_2(t)$ lie in the same connected component of
$X\setminus B_r(x)$ 
for all $t,r>0$ with $t\ge r$.
An equivalence class of cofinal rays is called an \textit{end} of $X$.
\end{defn}

\begin{proof}[Proof of Corollary~\ref{end}]
Suppose that the number of ends of $X$ is greater than $C^2+1$.
For any sequence $\epsilon_i\to 0$,
the space $(X,\epsilon_id,\mu)$ satisfies \textup{BG}$(0,n)$ with $C$
by Lemma~\ref{basiclem}.
For a point $x\in X$,
there exists a subsequence $\epsilon_j\to 0$ such that
$\{(X,\epsilon_jd,\mu,x)\}$ pointed measured Gromov--Hausdorff converges
to some pointed metric measure space $(X_{\infty},d_{\infty},\mu_{\infty},x_{\infty})$.
The limit space also satisfies \textup{BG}$(0,n)$ with $C$.
Since the number of ends of $X$ is greater than two,
the point $x_{\infty}$ is a local cut point.
Then the degree of $x_{\infty}$ is equal to the number of ends of $X$.
This contradicts Theorem~\ref{A}.
\end{proof}

\subsection{Branch points}

In this subsection,
we will give an obstruction condition for the existence of a local cut point
in a metric measure space $X$ satisfying BG$(k,n)$ with $1\le C<\sqrt{2}$.
Let us recall that the measured Gromov--Hausdorff limits of
manifolds $M_i$ with $\mathrm{Ric}_{M_i}\ge (n-1)k$ satisfy
BG$(k,n)$ with $C=1$.
It is conjectured that the limit space has no local cut point
unless the limit is one-dimensional.

Assume that $X$ has a local cut point $x$.
Let $x$ be an $r$-cut point.
Let $\gamma:[0,l]\to X$ be a minimal geodesic with $\gamma(l)=x$.
Assume that $l$ is sufficiently small ($l\le r/3$ for example).
If $\sigma:[0,L]\to X$ is a minimal geodesic such that
$\sigma(0)=\gamma(0)$ and $\sigma(L)\in B_l(\gamma(l))$ $(L\ge l)$,
then $\sigma(l)=\gamma(l)$, that is, $\sigma$ passes through $x$.
We then define two kinds of branch points of the geodesic $\gamma$.

\begin{defn}[branch point]\label{branchpt}\upshape
Let $\gamma$ be as above.
We say that $\gamma(l)$ is a \textit{branch point} of $\gamma$ if
for all $\epsilon>0$
there exist two distinct points $x_1, x_2\in B_{\epsilon}(\gamma(l))$
such that $d(\gamma(0),x_1)=d(\gamma(0),x_2)>l$.
\end{defn}

Note that $\gamma(l)$ is a branch point of $\gamma$
if and only if no neighborhood of $\gamma(l)$ is a segment.

\begin{figure}[tbp]
\begin{center}
 \includegraphics[width=10cm,clip]{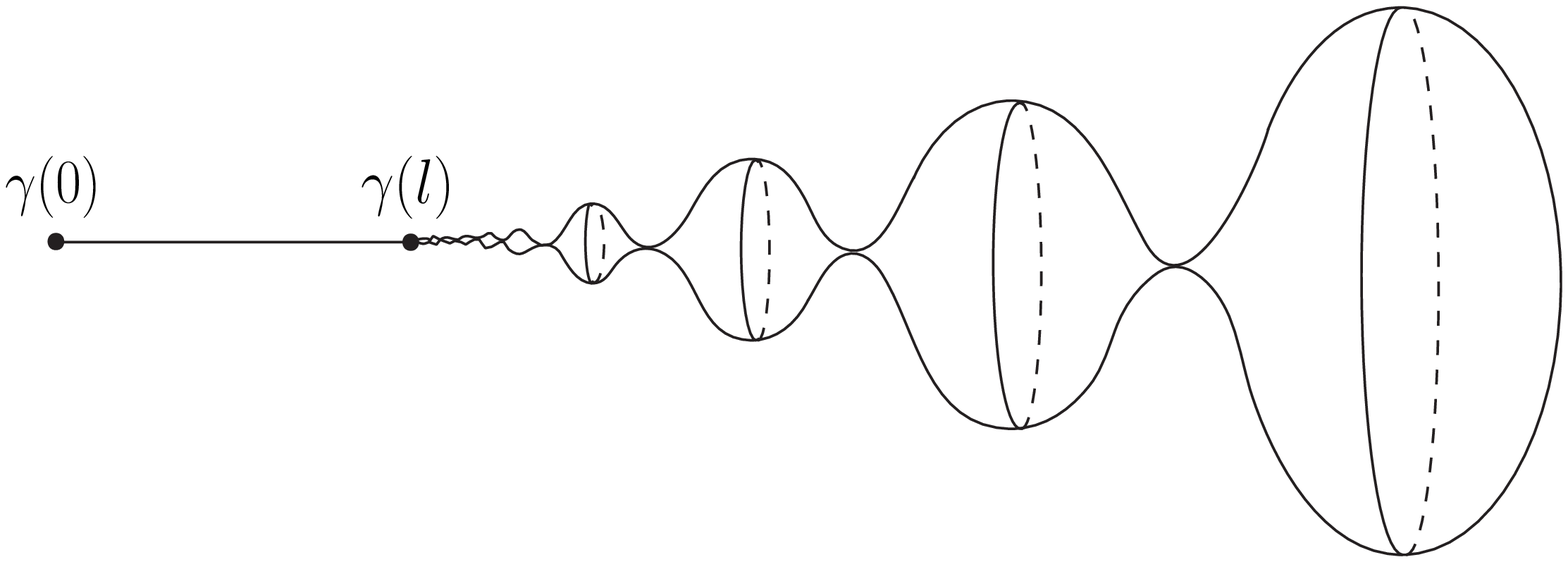}
\caption{Weak branch point}
\label{weakbranchpt}
\end{center}
\end{figure}

\begin{defn}[{weak branch point~\cite[Section 5]{CCII}}]\upshape
Let $\gamma$ be as in Definition~\ref{branchpt}.
We say that $\gamma(l)$ is a \textit{weak branch point} of $\gamma$ if
for all $\epsilon>0$
and all two points $x_1,x_2\in B_{\epsilon}(\gamma(l))$
with $d(\gamma(0),x_1), d(\gamma(0),x_2)>l$
there exist minimal geodesics $\sigma_i:[0,l_i]\to X$ $(i=1,2)$
from $\gamma(l)$ to $x_i$ such that
$\sigma_1(s)=\sigma_2(s)$ for some $s>0$.
See Figure~\ref{weakbranchpt}.
\end{defn}

\begin{ex}\upshape
(1)
Let $X=[a,b]$.
Every interior point $x\in (a,b)$ is a weak branch point 
of all paths $\gamma:[0,l]\to X$ with $\gamma(l)=x$ and $l\ll 1$.
The point $x$ is not a branch point of any $\gamma$.

(2)
Let $X$ be the length space in Figure~\ref{weakbranchpt}.
The intersection $x$ of two spheres is a branch point of all paths
$\gamma:[0,l]\to X$ with $\gamma(l)=x$ and $l\ll 1$.
However, the point $x$ is \textit{not} a weak branch point.
Although a branch point is not necessarily a weak branch point,
we use the term ``weak'' as in \cite[Section 5]{CCII}.
\end{ex}

Assume now that there exists a local cut point in 
a metric measure space satisfying BG$(k,n)$ with $1\le C<\sqrt{2}$.
We give an obstruction condition for the existence of a local cut point
by observing geodesics which pass through a local cut point.

\begin{thm}[cf.~{\cite[Theorem 5.1]{CCII}}]\label{weakbranch}
Let $(X,d,\mu)$ be a metric measure space satisfying 
the generalized Bishop--Gromov inequality 
\textup{BG}$(k,n)$ with constant $C$ 
for some $k\in \mathbb{R}$, $n\ge 1$, and $1\le C<\sqrt{2}$.
If there exists a local cut point $x$ in $X$,
then it is a weak branch point 
of all geodesics $\gamma:[0,l]\to X$ 
with $\gamma(l)=x$ and $l \ll 1$.
\end{thm}

The proof is essentially the same as in \cite{CCII}.
We give the proof for the completeness of the paper.

\begin{proof}
If $x$ has a neighborhood that is a segment,
then $x$ is a weak branch point of all paths which go to the point.
Assume now that no neighborhood of $x$ is a segment.
Let $\gamma:[0,l]\to X$ be a geodesic that branches at $\gamma(l)=x$,
where $l$ is sufficiently small.
Suppose that $x$ is not a weak branch point:
There exist two points $x_1,x_2\in B_l(\gamma(l))$
with $d(\gamma(0),x_1), d(\gamma(0),x_2)>l$ such that
for all geodesics $\sigma_i:[0,l_i]\to X$ $(i=1,2)$
from $\gamma(l)$ to $x_i$,
we have $\sigma_1(s)\neq \sigma_2(s)$ for all $0<s\le \min \{ l_1,l_2\}$.
We may assume that $d(\gamma(l),x_1)=d(\gamma(l),x_2)$,
and denote it by $l^{\prime}$.

Consider all minimal geodesics from $x$ to $x_i$ for $i=1,2$.
By the assumption, the union of all minimal geodesics from $x$ to $x_1$
and that of all minimal geodesics from $x$ to $x_2$ 
have no intersection except $x$.
We will use the generalized Bishop--Gromov inequality (\ref{BG(k,n)})
with the base point $x_i$.
Let $0<\epsilon \ll l^{\prime}$.
We denote by $U$ the connected component of 
$B_{\epsilon}(x)\setminus \{ x\}$ that contains $\gamma$.
We see that
$A_{l^{\prime},\,l^{\prime}+\epsilon}(x_i)$ contains $U$.
Set $S_i=S_{l^{\prime}-\epsilon,\,l^{\prime}}(x_i,U)$ for $i=1,2$.
Let us recall that $S_i$ is the intersection of $A_{l^{\prime}-\epsilon,\,l^{\prime}}(x_i)$
and all minimal geodesics from $x_i$ to each point in $U$.
Then $S_1\cap S_2$ is empty by the assumption.
For $i=1,2$ the generalized Bishop--Gromov inequality~(\ref{BG(k,n)}) induces
\begin{align}\label{S_i}
C^{-1} \frac{V_{k,\,n}(l^{\prime}-\epsilon,l^{\prime})}
     {V_{k,\,n}(l^{\prime},l^{\prime}+\epsilon)}
\le \frac{\mu\big(S_{l^{\prime}-\epsilon,\,l^{\prime}}(x_i,U)\big)}{\mu (U)}
= \frac{\mu(S_i)}{\mu(U)}.
\end{align}
Let $U^{\prime}$ denote the union of connected components of
$B_{\epsilon}(x)\setminus \{ x\}$ that contain $S_1$ and $S_2$.
We see that
$A_{l,\,l+\epsilon}(\gamma(0))$ contains $U^{\prime}$.
Since $S_i\subset U^{\prime}$ holds for $i=1,2$,
we have $\mu(S_1)+\mu(S_2)\le \mu(U^{\prime})$.
By summing up (\ref{S_i}) for $i=1,2$, it follows that
\begin{align}\label{S_1+S_2}
2C^{-1} \frac{V_{k,\,n}(l^{\prime}-\epsilon,l^{\prime})}
     {V_{k,\,n}(l^{\prime},l^{\prime}+\epsilon)}
\le \frac{\mu(S_1)+\mu(S_2)}{\mu (U)}
\le \frac{\mu(U^{\prime})}{\mu(U)}.
\end{align}

Applying (\ref{BG(k,n)}) 
with the base point $\gamma(0)$, we obtain
\begin{align}\label{Uprime}
\frac{\mu(U^{\prime})}
{\mu\big(S_{l-\epsilon,\,l}(\gamma(0),U^{\prime})\big)}
\le 
C\, \frac{V_{k,\,n}(l,l+\epsilon)}
     {V_{k,\,n}(l-\epsilon,l)}.
\end{align}
Since $U$ contains
$S_{l-\epsilon,\,l}(\gamma(0),U^{\prime})$,
it, together with (\ref{S_1+S_2}) and (\ref{Uprime}), follows that
\begin{align*}
2C^{-1} \frac{V_{k,\,n}(l^{\prime}-\epsilon,l^{\prime})}
       {V_{k,\,n}(l^{\prime},l^{\prime}+\epsilon)}
\le C\, \frac{V_{k,\,n}(l,l+\epsilon)}{V_{k,\,n}(l-\epsilon,l)}.
\end{align*}
Taking $\epsilon \to 0$, 
we obtain $2\le C^2$, which is a contradiction.
\end{proof}

\begin{figure}[tbp]
\begin{center}
 \includegraphics[width=9cm,clip]{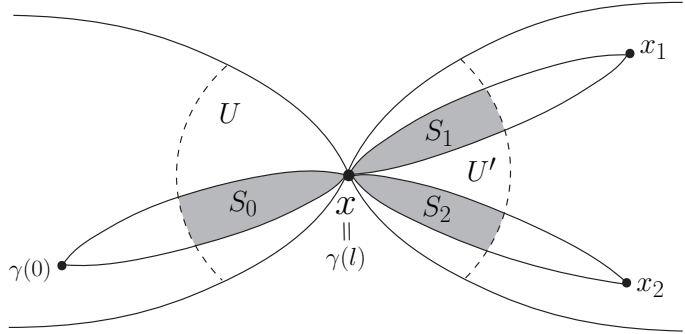}
\caption{Proof of Theorem~\ref{weakbranch}}
\label{Proofofweakbranch}
\end{center}
\end{figure}

\begin{rem}\upshape
In \cite[Theorem 5.1]{CCII},
Cheeger and Colding studied
the limit space which contains a one-dimensional piece
and which is not one-dimensional.
It follows from Theorem~\ref{weakbranch} that the length space in Figure~\ref{weakbranchpt}
can not satisfy BG$(k,n)$ with $1\le C<\sqrt{2}$ for any $k$, $n$, and any measure.
This does not follow from \cite[Theorem 5.1]{CCII}.
\end{rem}

As a corollary of Theorem~\ref{weakbranch}, 
we obtain the particular case of Theorem~\ref{A}.

\begin{thm}\label{partofThmA}
Let $(X,d,\mu)$ satisfy \textup{BG}$(k,n)$ with $C$
for some $k\in \mathbb{R}$, $n\ge 1$, and $1\le C<\sqrt{2}$.
Assume that there exists a local cut point $x$ in $X$.
Then we have $\deg(x)=2$.
\end{thm}

\begin{proof}
Suppose that $\deg(x)\ge 3$ holds.
For all geodesics $\gamma:[0,l]\to X$ with $\gamma(l)=x$ and $l\ll 1$,
the local cut point $x=\gamma(l)$ is not a weak branch point of $\gamma$.
This is because the geodesic $\gamma$ is extended to two connected components of
$B_r(x)\setminus \{ x\}$ which do not contain $\gamma$,
where $r$ is sufficiently small.
This contradicts Theorem~\ref{weakbranch}.
\end{proof}

\subsection{Local geometric structure}

Recall the definition of an $r$-cut point (Definition~\ref{cutpt}).
We investigate the geometric structure of the neighborhood of
an $r$-cut point in a metric measure space $X$ satisfying BG$(k,n)$ with $1\le C<\sqrt{2}$.

Assume that there exists an $r$-cut point $x$ in $X$.
We then have $\mathrm{diam}(O\cap S_r(x))\le \mathrm{diam}(\overline{B}_r(x))\le 2r$
for all connected components $O$ of $\overline{B}_r(x)\setminus \{ x \}$. 
By improving the method of the proof of Theorem~\ref{weakbranch},
we obtain a more precise estimate of the diameter as follows:

\begin{thm}\label{diam}
For all $k\in \mathbb{R}$, $n\ge 1$, $1\le C<\sqrt{2}$, and all $R>0$
there exists a constant $\delta=\delta(k,n,C,R)>0$ depending only on $k$, $n$, $C$, and $R$ 
such that the following holds:
Let $(X,d,\mu)$ be a metric measure space satisfying 
the generalized Bishop--Gromov inequality \textup{BG}$(k,n)$
with constant $C$.
If $X$ has an $r$-cut point $x$ with $0<r\le R$,
then 
$$\mathrm{diam}\big(O\cap S_r(x)\big)\le (2-\delta)r$$
holds for all connected components $O$ of $\overline{B}_r(x)\setminus \{ x \}$.

\end{thm}

\begin{figure}[tbp]
\begin{center}
 \includegraphics[width=7cm,clip]{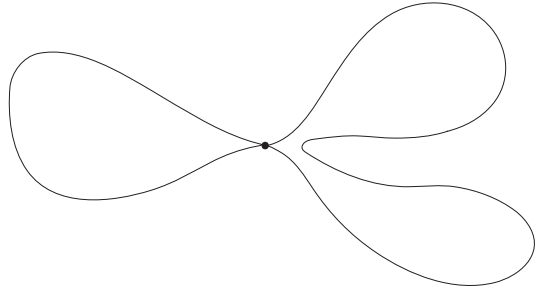}
\caption{``Three-pronged'' space}
\label{notBG}
\end{center}
\end{figure}

\begin{rem}\label{three-pronged}\upshape
The constant $\delta(k,n,C,R)$ in Theorem~\ref{diam} 
is independent of a metric measure space $(X,d,\mu)$.
Moreover, we can calculate the precise value;
see Remark~\ref{delta}.
It follows from Theorem~\ref{diam} that
for fixed $k,n$, and $C$
the metric space in Figure~\ref{notBG} can not satisfy 
BG$(k,n)$ with $C$ for any measure,
provided a narrow part is sufficiently small.
Roughly speaking, combined with Theorem~\ref{partofThmA},
``three-pronged'' spaces can not satisfy BG$(k,n)$ with $1\le C<\sqrt{2}$. 
\end{rem}

\begin{proof}[Proof of Theorem~\ref{diam}]
For $k\in \mathbb{R}$, $n\ge 1$, $1\le C<\sqrt{2}$, and $R>0$,
let $(X,d,\mu)$ be a metric measure space satisfying BG$(k,n)$ with $C$
and $x\in X$ an $r$-cut point with $0<r\le R$.
Fix $0<\delta \ll r$.
Suppose that
$$\mathrm{diam}\big(O\cap S_r(x)\big)>(2-\delta)r$$
holds for some connected component $O$ of $\overline{B}_r(x)\setminus \{ x \}$.
We take two points $x_1,x_2\in O\cap S_r(x)$ with $d(x_1,x_2)>(2-\delta)r$
and then choose minimal geodesics 
$\sigma_i:[0,r]\to X$ from $x$ to $x_i$ $(i=1,2)$.
Since $\deg(x)=2$,
we denote by $O^{\prime}$ another connected component of 
$\overline{B}_r(x)\setminus \{ x \}$.
For $0<\epsilon<\delta r$,
denote by $U$ the connected component of $B_{\epsilon}(x)\setminus \{ x \}$
which does not contain $\sigma_1$ and $\sigma_2$.
We see that
$A_{r/3,\,r/3+\epsilon}(\sigma_i(r/3))$ contains $U$
for $i=1,2$.
\begin{claim}
Each minimal geodesic from $\sigma_i(r/3)$ to each point in $U$ passes through $x$
$(i=1,2)$.
\end{claim}
\begin{proof}
Suppose that there exists a minimal geodesic $\gamma:[0,l]\to X$
from $\sigma_i(r/3)$ to some point $y$ in $U$
such that $\gamma$ does not pass through $x$.
Then, $\gamma(t)\notin B_r(x)$ for some $t\in [0,l]$
since $x$ is an $r$-cut point.
Therefore,
\begin{equation*}
d(y,\sigma_i(r/3))\ge d(y,\gamma(t))
\ge d(\gamma(t),x)-d(y,x)>r-\epsilon.
\end{equation*}
On the other hand, 
\begin{equation*}
d(y,\sigma_i(r/3))\le 
d(y,x)+d(x,\sigma_i(r/3))<\epsilon +\frac{r}{3}.
\end{equation*}
This is a contradiction.
\end{proof}   

Setting $S_i=S_{r/3-\delta r-\epsilon,\,r/3-\delta r}(\sigma_i(r/3),U)$ for $i=1,2$,
we have

\begin{claim}\label{empty}
$S_1\cap S_2=\emptyset$. 
\end{claim}
\begin{proof}
By the assumption, we have
\begin{align*}
(2-\delta)r<d(x_1,x_2)
&\le d(x_1,\sigma_1(r/3))+d(\sigma_1(r/3),\sigma_2(r/3))
+d(\sigma_2(r/3),x_2)\\
&=\frac{2r}{3}+d(\sigma_1(r/3),\sigma_2(r/3))+\frac{2r}{3}.
\end{align*}
It follows that
$(2/3-\delta) r<d(\sigma_1(r/3),\sigma_2(r/3)).$
Suppose that there exists a point $y\in S_1\cap S_2$. 
We then have
\begin{equation*}
d(\sigma_1(r/3),\sigma_2(r/3))
\le d(\sigma_1(r/3),y)+d(y,\sigma_2(r/3))
\le 2\Big(\frac{r}{3}-\delta r\Big) =\Big(\frac{2}{3}-2\delta \Big) r.
\end{equation*}
This is a contradiction.
\end{proof}

Using the generalized Bishop--Gromov inequality~(\ref{BG(k,n)}) 
with the base point $\sigma_i(r/3)$,
we obtain for $i=1,2$
\begin{align}\label{BG-S_i}
C^{-1}\frac{V_{k,\,n}\big((1/3-\delta)r-\epsilon,(1/3-\delta)r\big)}
{V_{k,\,n}(r/3,r/3+\epsilon)}
\le \frac{\mu(S_i)}{\mu(U)}.
\end{align}
We set $U^{\prime}=S_1\cup S_2$.
Since $S_1\cap S_2=\emptyset$ (Claim~\ref{empty}),
we have $\mu(S_1)+\mu(S_2)\le \mu(U^{\prime})$.
Therefore, summing up (\ref{BG-S_i}) for $i=1,2$ gives
\begin{align}\label{S_1+S_2diam}
2C^{-1}\frac{V_{k,\,n}\big((1/3-\delta)r-\epsilon,(1/3-\delta)r\big)}
{V_{k,\,n}(r/3,r/3+\epsilon)}
\le \frac{\mu(S_1)+\mu(S_2)}{\mu(U)}
\le \frac{\mu(U^{\prime})}{\mu(U)}.
\end{align}

Next, we take an arbitrary point $x_0$ in $O^{\prime}\cap S_{r/3}(x)$.
Such a point exists since $x$ is an $r$-cut point.

\begin{claim}
Each minimal geodesic from $x_0$ to each point in $U^{\prime}$ passes through $x$.
\end{claim}
\begin{proof}
Suppose that there exists a minimal geodesic $\gamma:[0,l]\to X$
from $x_0$ to some point $y$ in $U^{\prime}$ such that $\gamma$ does not pass through $x$.
Then, $\gamma(t)\notin B_r(x)$ for some $t\in [0,l]$.
Therefore,
\begin{equation*}
d(x_0,y)\ge d(x_0,\gamma(t_0))\ge d(\gamma(t_0),x)-d(x,x_0)
\ge r-\frac{r}{3}=\frac{2r}{3}.
\end{equation*}
On the other hand,
\begin{equation*}
d(x_0,y)\le d(x_0,x)+d(x,y)\le \frac{r}{3}+\delta r+\epsilon.
\end{equation*}
This is a contradiction.
\end{proof}

Set $S_0=S_{r/3-\epsilon,\,r/3}(x_0,U^{\prime})$.
The generalized Bishop--Gromov inequality~(\ref{BG(k,n)}) 
with the base point $x_0$ implies
\begin{align}\label{Uprimediam}
\frac{\mu(U^{\prime})}{\mu(S_0)}\le
C\, \frac{V_{k,\,n}\big((1/3+\delta)r,(1/3+\delta)r+\epsilon\big)}
{V_{k,\,n}(r/3-\epsilon,r/3)}.
\end{align}
Since $U$ contains $S_0$,
inequalities (\ref{S_1+S_2diam}) and (\ref{Uprimediam}) give
\begin{align}\label{lasteq}
2\le C^2\, \frac{V_{k,\,n}(r/3,r/3+\epsilon)}
{V_{k,\,n}(r/3-\epsilon,r/3)}
\frac{V_{k,\,n}\big((r/3+\delta)r,(r/3+\delta)r+\epsilon\big)}
{V_{k,\,n}\big((r/3-\delta)r-\epsilon,(r/3-\delta)r\big)}.
\end{align}
After letting $\epsilon \to 0$,
a sufficiently small $\delta>0$ implies a contradiction.
This completes the proof of Theorem~\ref{diam}.
\end{proof}

\begin{figure}[tbp]
\begin{center}
 \includegraphics[width=10cm,clip]{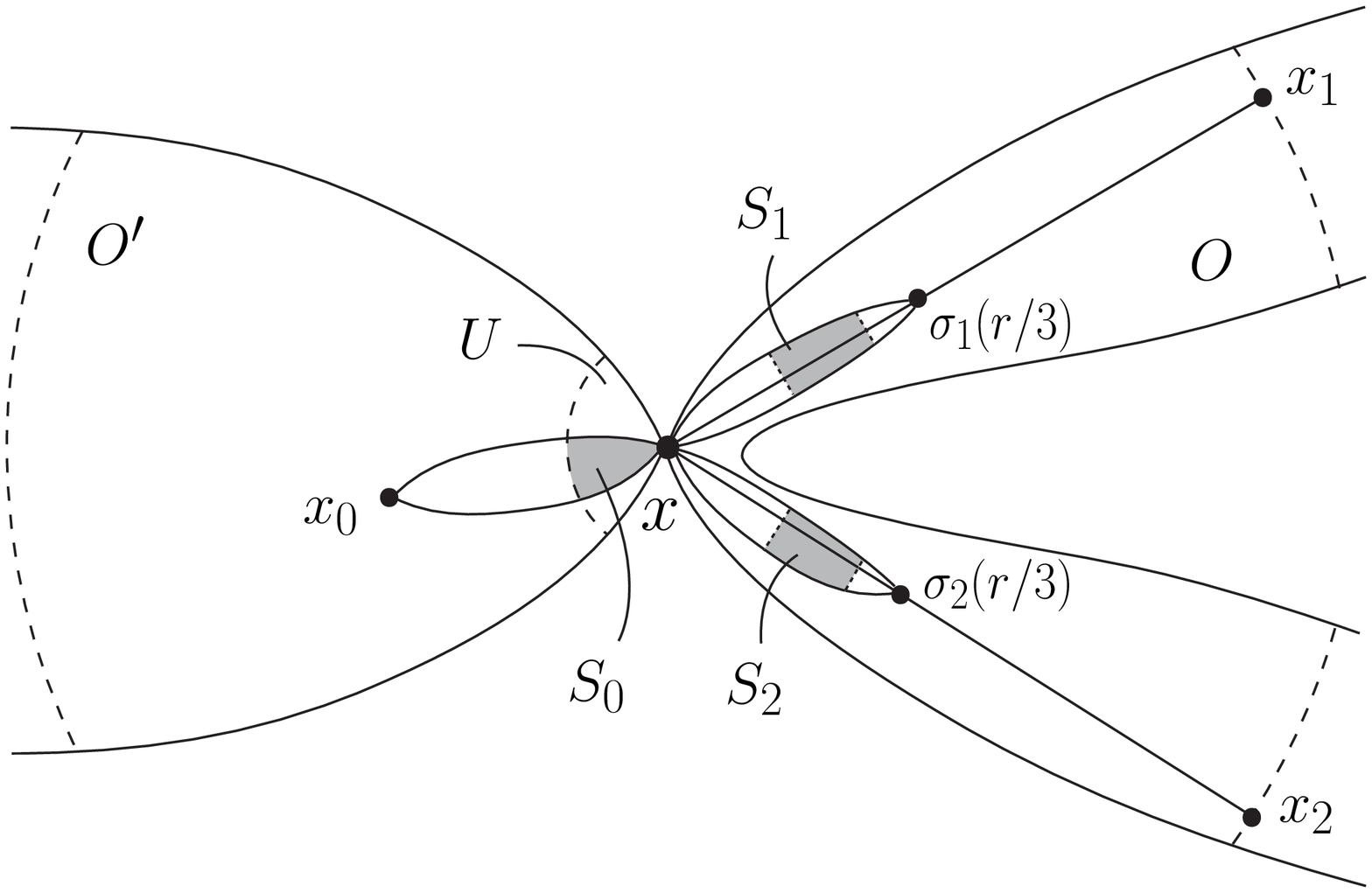}
\caption{Proof of Theorem~\ref{diam}}
\label{Proofdiam}
\end{center}
\end{figure}

\begin{rem}\label{delta}\upshape
In the case $k=0$, we have
$$\delta(0,n,C,R)=\frac{1}{4}\ \frac{(2/C^2)^{1/(n-1)}-1}{(2/C^2)^{1/(n-1)}+1}.$$
Indeed, the right-hand side of (\ref{lasteq}) is equal to
$$C^2\, \frac{(r/3+\epsilon)^n-(r/3)^n}{(r/3)^n-(r/3-\epsilon)^n}\,
\frac{[ (1/3+\delta)r+\epsilon]^n-[ (1/3+\delta)r]^n}
{[ (1/3-\delta)r]^n-[ (1/3-\delta)r-\epsilon]^n},$$
which converges to $C^2[ (1/3+\delta)/(1/3-\delta)]^{n-1}$
as $\epsilon \to 0$.
Therefore, it suffices to determine a positive number $\delta$
such that $C^2[ (1/3+\delta)/(1/3-\delta)]^{n-1}<2$.
\end{rem}

\subsection{Convergence of local cut points}

We observe the structure of the accumulation of local cut points 
in a metric measure space satisfying
BG$(k,n)$ with $1\le C<\sqrt{2}$
by using Theorem~\ref{diam}.

Let $(X,d,\mu)$ satisfy \textup{BG}$(k,n)$ with $C$
for some $k\in \mathbb{R}$, $n\ge 1$, and $1\le C<\sqrt{2}$.
Assume that there exist three $r$-cut points $x_1,x_2,x_3\in X$ ($r>0$)
such that $d(x_i,x_j)$ is sufficiently small for all $1\le i,j\le 3$.
We recall that $\deg(x_i)=2$ (Theorem~\ref{partofThmA}).
Denote by $O_1,O_1^{\prime}$ 
the connected components of $\overline{B}_r(x_1)\setminus \{ x_1\}$.
Assume then that $O_1^{\prime}$ contains $x_2$ and $x_3$
and we have $d(x_1,x_2)<d(x_1,x_3)$ without loss of generality.
For $i=2,3$
denote by $O_i$ the connected component of
$\overline{B}_r(x_i)\setminus \{ x_i\}$ which contains $x_1$
and by $O_i^{\prime}$ another one.
By the definition,
$O_2\cap O_3$ is nonempty.
We say that $\{ x_1,x_2,x_3\}$ \textit{stands in a line} if
$O_2^{\prime}\cap O_3^{\prime}$ is nonempty.

As a corollary of Theorem~\ref{diam}, we have the following:

\begin{cor}\label{standline}
Let $(X,d,\mu)$ be a metric measure space satisfying 
\textup{BG}$(k,n)$ with $C$ for some $k\in \mathbb{R}$, $n\ge 1$,
and $1\le C<\sqrt{2}$.
Let $\delta=\delta(k,n,C,R)$ be the constant in Theorem~\ref{diam} for $R>0$.
Assume that there exist $r$-cut points $x_1,x_2,x_3\in X$ with $0<r\le R$.
If $d(x_i,x_j)<\delta r/6$ holds for all $1\le i,j\le 3$,
then $\{ x_1,x_2,x_3\}$ stands in a line.
\end{cor}

\begin{figure}[tbp]
\begin{center}
 \includegraphics[width=10cm,clip]{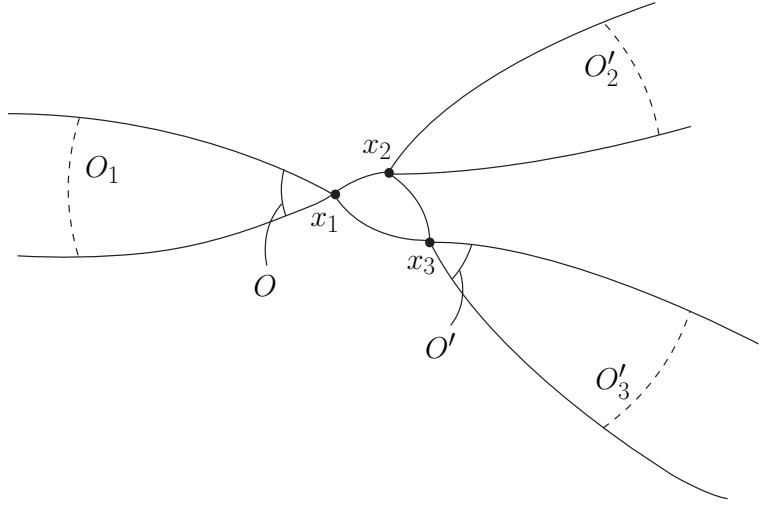}
\caption{Proof of Corollary~\ref{standline}}
\label{Proofsequence}
\end{center}
\end{figure}

\begin{proof}
Suppose that $O_2^{\prime}\cap O_3^{\prime}$ is empty;
see Figure~\ref{Proofsequence}.
Since $x_2$ is an $r/3$-cut point, there exist at least two connected components of
$O_2\cap S_{r/3}(x_2)=(O_1\cap S_{r/3}(x_2))\cup (O_3^{\prime}\cap S_{r/3}(x_2))$.
Denote by $O$ one of connected components of $O_1\cap S_{r/3}(x_2)$
and by $O^{\prime}$ one of connected components of $O_3^{\prime}\cap S_{r/3}(x_2)$.

Since $x_2$ is an $r$-cut point,
all minimal geodesics from each point in $O$ to each point in  $O^{\prime}$
pass through $x_1$ and $x_3$.
Therefore, we have
\begin{align}\label{OOprime}
d(O,O^{\prime})=d(x_1,O)+d(x_1,x_3)+d(x_3,O^{\prime}).
\end{align}
Since $d(x_i,x_j)<\delta r/6$ holds by the assumption, we have
\begin{align}
d(x_1,O)\ge \frac{r}{3}-d(x_2,x_1)>\frac{r}{3}-\frac{\delta r}{6}
        =\bigg( 1-\frac{\delta}{2}\bigg) \frac{r}{3}\label{x_1}. 
\end{align}
Similarly, we have
\begin{align}
d(x_3,O^{\prime})>\bigg( 1-\frac{\delta}{2}\bigg) \frac{r}{3}\label{x_3}. 
\end{align}
Hence, relations (\ref{OOprime}), (\ref{x_1}) and (\ref{x_3}) imply
\begin{align*}
\mathrm{diam}(O_2\cap S_{r/3}(x_2))
&\ge d(O,O^{\prime})\\
&=d(x_1,O)+d(x_1,x_3)+d(x_3,O^{\prime})\\
&>\bigg( 1-\frac{\delta}{2} \bigg)\frac{r}{3}+\bigg( 1-\frac{\delta}{2} \bigg)\frac{r}{3}\\
&=\big( 2-\delta \big) \frac{r}{3}.
\end{align*}
This contradicts Theorem~\ref{diam}.
\end{proof}

Assume now that 
there exists a sequence $\{ x_i\}_{i=1}^{\infty}$
of $r$-cut points in $X$
such that $d(x_i,x_j)$ is sufficiently small for all $i$ and $j$.
Denote by $O_i,O_i^{\prime}$
the connected components of $\overline{B}_r(x_i)\setminus \{ x_i\}$.
Assume that $O_i$ contains $x_j$ for all $j<i$.
We say that $\{ x_i\}_{i=1}^{\infty}$ \textit{stands in a line} if
$O_i^{\prime}\cap O_{i+1}^{\prime}$ is nonempty for all $i$.

Corollary~\ref{standline} implies

\begin{cor}\label{closed}
Let $(X,d,\mu)$ be as in Corollary~\ref{standline}.
For each $r>0$,
the set of all $r$-cut points in $X$ is closed.
\end{cor}

\begin{proof}
Assume that there exist a sequence $\{x_i\}_{i=1}^{\infty}$ of $r$-cut points in $X$, 
which converges to a point $x$ in $X$.
By Corollary~\ref{standline}, 
$\{ x_i\}_{i=N}^{\infty}$ stands in a line for a sufficiently large $N$.
Hence the limit $x$ is a local cut point.
Since each $x_i$ is an $r$-cut point,
so is $x$.
\end{proof}

\begin{rem}\upshape
Let $\{ x_i\}_{i=1}^{\infty}$ be a sequence of $r_i$-cut points,
where $r_i\to 0$.
It is then not necessarily that a limit point of $\{ x_i\}$ is a local cut point. 
Consider the sequence $\{ 1/i\}_{i=1}^{\infty}\subset [0,1]$.
\end{rem}

In fact,
Corollary~\ref{closed} holds without the assumption 
BG$(k,n)$ with $1\le C<\sqrt{2}$.

\begin{prop}\label{loccpt}
Let $(X,d)$ be a complete, locally compact length space.
Assume that for $r>0$ there exists a sequence $\{ x_i\}_{i=1}^{\infty}$ of 
$r$-cut points in $X$ 
which converges to a point $x\in X$.
Then,
for some sufficiently large $N$,
we have $\deg(x_i)=2$ for all $i\ge N$. 
Moreover,
$\{ x_i\}_{i=N}^{\infty}$ stands in a line.
$($Hence, the limit point $x$ is also an $r$-cut point.
We do not always have $\deg(x)=2$.$)$
\end{prop}

\begin{proof}
Suppose that for each $N$ there exists $i\ge N$
such that $\deg(x_i)\ge 3$.
We may assume that $\deg(x_i)\ge 3$ for all $i$,
and $d(x,x_i)>d(x,x_j)$ for all $i<j$.
For each $i$,
there exists a connected component $O_i$ of $\overline{B}_r(x_i)\setminus \{ x_i\}$
which does not contain $x$ and $x_{i-1}$.
Since $O_i\cap S_r(x_i)$ is nonempty,
we have a contradiction to the locally compactness at $x$.
Thus, it follows that $\deg(x_i)=2$ for all sufficiently large $i$.

Next, we show that 
$\{ x_i\}_{i=N}^{\infty}$ stands in a line 
for some sufficiently large $N$.
Suppose that $\{ x_i\}_{i=N}^{\infty}$ does not stand in a line
for any large $N$.
By taking a subsequence,
we may assume that 
$O_i^{\prime}\cap O_{i+1}^{\prime}=\emptyset$
for $i=1,3,5,\dots$.
Since $O_i^{\prime}\cap S_r(x_i)$ is nonempty,
we have a contradiction to the locally compactness at $x$. 
\end{proof}

In general,
Proposition~\ref{loccpt} does not necessarily hold 
without the locally compactness of $X$
\begin{ex}\upshape
Consider the set
$\bigcup_{i=1}^{\infty}\{ (2/i,y)\,|\, 0\le y\le 2\}
\cup \{ (x,0)\,|\, 0\le x\le 2\}\subset \mathbb{R}^2$
with the induced distance.
Although the point $(2/i,1/i)$ is a $1$-cut point for every $i\in \mathbb{N}$,
the limit point $(0,0)$ is not a local cut point.
\end{ex}

\section{The Poincar\'{e} inequality}

In this section we prove Theorem~\ref{B}.
We begin by recalling the definition of a Poincar\'{e} inequality of type $(1,p)$.
Let $(X,d)$ be a metric space.

\begin{defn}[upper gradient~\cite{HeiK}]\upshape
Let $u$ be a function on $X$.
A Borel function $g:X\to [0,\infty]$ is called an \textit{upper gradient} of $u$
if for all paths $\gamma:[0,l]\to X$ (proportional to arclength), we have
$$|u(\gamma(l))-u(\gamma(0))|\le \int_0^lg(\gamma(t))\ dt.$$
\end{defn}

Every function has an upper gradient $g\equiv \infty$,
and hence upper gradients are never unique.
For a Lipschitz function $u:X\to \mathbb{R}$,
we define $|\nabla u|:X\to \mathbb{R}$ by
$$|\nabla u|(x)=\limsup_{y\to x}\frac{|u(x)-u(y)|}{d(x,y)}$$
if $x$ is not isolated, and $|\nabla u|(x)=0$ if $x$ is isolated. 

\begin{prop}\label{upper}
If $u$ is a Lipschitz function,
then $|\nabla u|$ is an upper gradient of $u$.
\end{prop}

\begin{proof}
Let $\gamma:[0,l]\to X$ be a path.
The function $u\circ \gamma$ is Lipschitz
and hence is differentiable almost everywhere by Rademacher's theorem.
Since $\gamma$ is proportional to arclength,
$$|(u\circ \gamma)^{\prime}(t)|
= \lim_{s\to t}\frac{|u(\gamma(t))-u(\gamma(s))|}{|t-s|}
\le \lim_{s\to t}\frac{|u(\gamma(t))-u(\gamma(s))|}{d(\gamma(t),\gamma(s))}
\le |\nabla u|(\gamma(t))$$
holds for all differentiable points $t\in [0,l]$.
Therefore
$$|u(\gamma(l))-u(\gamma(0))|\le \int_0^l|(u\circ \gamma)^{\prime}(t)|\ dt
\le \int_0^l|\nabla u|(\gamma(t))\ dt,$$
which completes the proof.
\end{proof}

Assume that $(X,d,\mu)$ is a complete, locally compact length space
equipped with a Borel measure such that 
$0<\mu(B_r(x))<+\infty$ holds for all $x\in X$ and all $0<r<+\infty$.
We denote
\begin{equation*}
u_B=\mint_Bu\ d\mu=\frac{1}{\mu(B)}\int_Bu\ d\mu
\end{equation*}
for $B\subset X$.

\begin{defn}[Poincar\'{e} inequality]\upshape
Let $1\le p<\infty$.
We say that $(X,d,\mu)$ satisfies
a \textit{Poincar\'{e} inequality of type} $(1,p)$,
if for all $R>0$ there exists a constant $C_P=C_P(p,R)>0$ 
depending only on $p$ and $R$ such that
\begin{align*}
\mint_{B_r(x)} | u-u_{B_r(x)}|\ d\mu 
\le C_P\,r
\bigg(\mint_{B_r(x)} g^p\ d\mu\bigg)^{1/p}
\end{align*}
holds for all $x\in X$, all $0<r\le R$, 
all measurable functions $u$, and all upper gradients $g$ of $u$.
\end{defn}

In our setting ($X$ is a length space),
a Poincar\'{e} inequality of type $(1,p)$ is derived from
a ``weak'' Poincar\'{e} inequality of type $(1,p)$
if we assume that $\mu$ is doubling
(see \cite{HajK} for details).

\begin{rem}\upshape
It follows from H\"{o}lder's inequality that
each metric measure space satisfying a Poincar\'{e} inequality of type $(1,p)$
also satisfies that of type $(1,q)$ 
for all $q\ge p$.

Keith and Zhong~\cite{KZ} proved the following:
Let $p>1$.
If $\mu$ is doubling and if $(X,d,\mu)$ satisfies a Poincar\'{e} inequality of type $(1,p)$,
then there exists $\epsilon>0$ such that $(X,d,\mu)$ satisfies
a Poincar\'{e} inequality of type $(1,q)$ for all $q>p-\epsilon$.
\end{rem}

\begin{proof}[Proof of Theorem~\ref{B}]
The proof is by contradiction;
suppose that there exists a local cut point $x_0$ 
satisfying (\ref{katei}).
Fix a sufficiently small $r>0$ such that 
$B_r(x_0)\setminus \{ x_0 \}$ is disconnected.
We choose two connected components $O_1,O_2$ of $B_r(x_0)\setminus \{ x_0 \}$.
For sufficiently large numbers $N\in \mathbb{N}$,
we define functions $u_N:B_r(x_0)\to \mathbb{R}$ as follows.
We set $U_i=O_i\cap (B_r(x_0)\setminus B_{1/N\mu(U_i)}(x_0))$ 
for $i=1,2$,
and define
\begin{equation*}
u_N(x)=\begin{cases}
            (-1)^{i+1}/\mu(O_i)  & \mbox{on}\ U_i,\\
            (-1)^{i+1}Nd(x_0,x)  & \mbox{on}\ O_i\setminus U_i,\\
            0                    & \mbox{on}\ B_r(x_0)\setminus (O_1\cup O_2). 
            \end{cases}
\end{equation*} 
The function $u_N$ is Lipschitz; hence,
$|\nabla u_N|$ is an upper gradient of $u_N$ (Proposition~\ref{upper}).
Since $(X,\mu)$ satisfies a Poincar\'{e} inequality of type $(1,p)$,
for $R\ge r$ there exists a constant $C_P=C_P(p,R)>0$ such that
\begin{align}\label{Poincare-u_N}
\mint_{B_r(x_0)} \big| u_N-(u_N)_{B_r(x_0)}\big|\ d\mu 
\le C_P\,r
\bigg(\mint_{B_r(x_0)} |\nabla u_N|^p\ d\mu\bigg)^{1/p}.
\end{align}

We first estimate the left-hand side of (\ref{Poincare-u_N}) from below.
We have
\begin{align*}
\hspace{5pt}
\mint_{B_r(x_0)} \big| u_N-(u_N)_{B_r(x_0)}\big| \ d\mu \phantom{\int}\\
&\hspace{-13em}\ge \mint_{B_r(x_0)}|u_N|\ d\mu 
 -\big| (u_N)_{B_r(x_0)}\big| \phantom{\int}\\
&\hspace{-13em}=
 \frac{1}{\mu(B_r(x_0))} 
 \Bigg[ \bigg( \frac{\mu(U_1)}{\mu(O_1)}+\int_{O_1\setminus U_1}Nd(x_0,x)\ d\mu(x) 
 +\frac{\mu(U_2)}{\mu(O_2)}+\int_{O_2\setminus U_2}Nd(x_0,x)\ d\mu(x) \bigg)\\
&\hspace{-9em}-\bigg| 
 \frac{\mu(U_1)}{\mu(O_1)}+\int_{O_1\setminus U_1}Nd(x_0,x)\ d\mu(x) 
 -\frac{\mu(U_2)}{\mu(O_2)}-\int_{O_2\setminus U_2}Nd(x_0,x)\ d\mu(x) \bigg| \, \Bigg] \\
&\hspace{-13em}\ge \frac{1}{\mu(B_r(x_0))} 
 \Bigg[ \bigg( \frac{\mu(U_1)}{\mu(O_1)}+\int_{O_1\setminus U_1}Nd(x_0,x)\ d\mu(x) 
 +\frac{\mu(U_2)}{\mu(O_2)}+\int_{O_2\setminus U_2}Nd(x_0,x)\ d\mu(x) \bigg) \\
&\hspace{-9em}-
 \bigg( \bigg| \frac{\mu(U_1)}{\mu(O_1)}-\frac{\mu(U_2)}{\mu(O_2)}\bigg|
 +\int_{O_1\setminus U_1}Nd(x_0,x)\ d\mu(x) 
 +\int_{O_2\setminus U_2}Nd(x_0,x)\ d\mu(x) \bigg) \Bigg] \\
&\hspace{-13em}= \frac{1}{\mu(B_r(x_0))}
 \bigg( \frac{\mu(U_1)}{\mu(O_1)}+\frac{\mu(U_2)}{\mu(O_2)}
 -\bigg| \frac{\mu(U_1)}{\mu(O_1)}-\frac{\mu(U_2)}{\mu(O_2)}\bigg| \bigg).
\end{align*}
Since $N$ is sufficiently large,
$\mu(U_i)/\mu(O_i)$ is approximately equal to one ($i=1,2$). 
Therefore, 
we see that the left-hand side of (\ref{Poincare-u_N}) is bounded below by 
the positive number $\mu(B_r(x_0))^{-1}$
which is independent of $N$.

Next, we estimate the right-hand side of (\ref{Poincare-u_N}) from above.
Note that
\begin{equation*}
|\nabla u_N|=\begin{cases}
              N & \mbox{on}\ (O_1\setminus U_1)\cup (O_2\setminus U_2),\\
              0 & \mbox{otherwise},
             \end{cases}
\end{equation*} 
and
$(O_1\setminus U_1)\cup (O_2\setminus U_2)$ is contained in
$B_{\max \{ 1/N\mu(O_1),1/N\mu(O_2) \}}(x_0)$.
It follows that
\begin{align*}
C_P\,r\bigg(\mint_{B_r(x_0)} |\nabla u_N|^p\ d\mu\bigg)^{1/p}
&= C_P\,r\Big[ \mu(B_r(x_0))^{-1}
 N^p\mu \big((O_1\setminus U_1)\cup (O_2\setminus U_2) \big) \Big]^{1/p} \\
&\le C_P\,r\Big[ \mu(B_r(x_0))^{-1}
 N^p\mu \big( B_{\max \{ 1/N\mu(O_1),1/N\mu(O_2) \}}(x_0)\big) \Big]^{1/p},
\end{align*}
which converges to zero as $N\to \infty$,
by the assumption (\ref{katei}).
This is a contradiction.
\end{proof}

Let $\alpha>0$.
Consider the metric measure space $\{ (x_1,x_2,\dots,x_n)\in \mathbb{R}^n\, |
\, x_1^2+x_2^2+\cdots +x_{n-1}^2\le x_n^{2\alpha}\}$
equipped with the Euclidean distance and the $n$-dimensional Lebesgue measure $\mathcal{L}^n$.
The origin $o=(0,0,\dots,0)$ is a local cut point.
We have $\mathcal{L}^n(B_r(o))=c(n)\,r^{\alpha(n-1)+1}$.
It follows from Theorem~\ref{B} that,
if the metric measure space satisfies
a Poincar\'{e} inequality of type $(1,p)$, then we have
$$\alpha\le \frac{p-1}{n-1}.$$


\end{document}